\documentclass[11pt]{article}
\usepackage[utf8]{inputenc}
\usepackage[margin=1in,left=1in]{geometry}
\usepackage{setspace}
\usepackage{graphicx}
\usepackage{placeins}   
\usepackage{float}
\usepackage{amssymb, amsmath, amsthm, graphics}
\usepackage{mathabx,epsfig}
\usepackage[mathscr]{euscript}
\usepackage{tikz}
\usepackage{pigpen}
\usepackage{bm}
\usetikzlibrary{calc}
\usetikzlibrary{matrix}
\usepackage{latexsym}
\usepackage{color}

\usepackage{enumitem}

\usepackage{tikz}
\input xy
\xyoption{all}
\usepackage{euscript}
\usepackage{multirow}
\usepackage{etex, pictexwd,dcpic}
\usetikzlibrary{positioning}
\usetikzlibrary{shapes.geometric}
\usetikzlibrary{shapes.misc}
\usetikzlibrary{calc}
\usetikzlibrary{positioning}

 \usepackage[arc,all]{xy}

\usepackage{pgflibraryarrows}
\usepackage{pgflibrarysnakes}

\usetikzlibrary{trees} 
\usetikzlibrary[trees] 

\newtheorem{theorem}{Theorem}
\newtheorem{definition}[theorem]{Definition}
\newtheorem{lemma}[theorem]{Lemma}
\newtheorem{corollary}[theorem]{Corollary}
\newtheorem{proposition}[theorem]{Proposition}

\newtheorem{remark}{Remark}

\newtheorem{example}{Example}

\numberwithin{equation}{section}

\renewcommand{\(}{\begin{equation*}}
\renewcommand{\)}{\end{equation*}}
\newcommand{\bea}{\begin{eqnarray*}}
\newcommand{\eea}{\end{eqnarray*}}

\renewcommand{\a}{\alpha}
\renewcommand{\b}{\beta}

\def\endofproof {\hfill{$\Box$}\\}

\newcommand{\cL}{\ensuremath{\mathcal L}}

\def\L{\ensuremath{{\cal L}}}

\newcommand{\Exterior}{\scalebox{.8}{\ensuremath \bigwedge}}

\def\L{\ensuremath{{\cal L}}}

\newcommand{\beq}{\begin{equation}}
\newcommand{\eeq}{\end{equation}}

\newcommand{\onto}{\twoheadrightarrow}

\newcommand{\into}{\hookrightarrow}

\newcommand{\op}[1]{\ensuremath{\operatorname{#1}}}

\newcommand{\theproof}{\noindent {\bf Proof.\ }}

\numberwithin{equation}{section}

\renewcommand{\(}{\begin{equation}}
\renewcommand{\)}{\end{equation}}

\def\1{{\bf 1}}

\def\<{\langle}
\def\>{\rangle}

\def\a{\alpha}
\def\b{\beta}

\numberwithin{equation}{section}


\newcommand{\RR}{\ensuremath{\mathbb R}}

\newcommand{\ZZ}{\ensuremath{\mathbb Z}}
\newcommand{\Z}{\ensuremath{\mathbb Z}} 
\newcommand{\QQ}{\ensuremath{\mathbb Q}}

\newcommand{\BB}{\ensuremath{\mathbf B}}

\newcommand{\CC}{\ensuremath{\mathbb C}}

\newcommand{\chp}{\ensuremath{\mathscr{C}\mathrm{h}^{+}}}

\newcommand{\Sp}{\ensuremath{\mathscr{S}\mathrm{p}}}
\newcommand{\sset}{\ensuremath{\mathrm{s}\mathscr{S}\mathrm{et}}}

\newcommand{\mf}{\ensuremath{\mathscr{M}\mathrm{f}}}

\newcommand{\sh}{\ensuremath{\mathscr{S}\mathrm{h}}}

\newcommand{\C}{\ensuremath{\mathscr{C}}}

\newcommand{\E}{\ensuremath{\mathscr{E}}}
\newcommand{\R}{\ensuremath{\mathscr{R}}}

\newcommand{\map}{\mathrm{Map}}

\begin{document}

\title{
Twisted differential  generalized cohomology theories \\
 and their Atiyah-Hirzebruch spectral sequence}

 \author{
 Daniel Grady and Hisham Sati\\
  }

\maketitle

\begin{abstract} 
We construct the Atiyah-Hirzebruch spectral sequence (AHSS) for 
twisted differential  generalized cohomology theories. This generalizes
to the twisted setting the authors' corresponding earlier construction for differential 
cohomology theories, as well as to the differential setting the AHSS 
for twisted generalized cohomology theories, including that of twisted K-theory 
by Rosenberg and Atiyah-Segal. In describing twisted  differential spectra 
we build on the work of Bunke-Nikolaus, but we find it useful for our purposes 
to take an approach that highlights direct analogies with classical bundles and 
that is at the same time amenable for calculations. We will, in particular, 
establish that twisted differential spectra are bundles of spectra equipped 
with a flat connection. Our prominent case will be twisted differential K-theory, 
for which we work out the differentials in detail. This involves differential 
refinements of primary and secondary cohomology operations the authors 
developed in earlier papers. We illustrate our constructions and computational 
tools with examples.
 \end{abstract}

 \tableofcontents
 
\section{Introduction}

There has been a lot of research activity in construction and using twisted cohomology theories.  The main example has been twisted K-theory $K^*(X; h)$ of a space $X$, where the twist $H$ is a cohomology class $h \in H^3(X, \Z)$ (see 
\cite{DK} \cite{Ro} \cite{BCMMS} 
\cite{AR}  \cite{TXL}
\cite{AS1} \cite{AS2} \cite{FHT} \cite{Ku}
\cite{ABG} \cite{AGG}   \cite{Ka2}   \cite{Go}).  
Periodic de Rham cohomology may be twisted by any odd degree cohomology class \cite{teleman} \cite{BSS} \cite{MW}. Twisted de Rham cohomology has roots that 
go back to \cite{RW}, and has attracted attention more recently -- see e.g. 
\cite{MW} \cite{S1} \cite{S2}. TMF can be twisted by a degree four class \cite{ABG} 
and Morava K-theory and E-theory admit twists by cohomology classes for 
every prime $p$ and any chromatic level $n$ \cite{SW}. 
More exotic spectra can be twisted in an even more unexpected way, 
including iterated algebraic K-theory of the connective K-theory spectrum \cite{LSW}. 

\medskip
Twisted generalized cohomology theories will admit an isomorphism with
 twisted periodic de Rham cohomology under the generalized twisted Chern 
 character (or Chern-Dold character). For degree-three 
 twisted K-theory,  the twisted Chern character appears from 
 various points of view in 
  \cite{AS2}\cite{Ka1} \cite{MS1} \cite{MS2}  \cite{TX} \cite{BGNT} \cite{Ka2} \cite{GT} \cite{HM}.  
  Chern characters can be considered for higher twisted theories, such 
 as Morava K-theory and E-theory \cite{SW} and iterated algebraic K-theory in \cite{LSW}.   Differential refinements for a twisted generalized cohomology theory 
are considered in \cite{BN}. 

\medskip
On the other hand, differential cohomology has become an active area of research 
(e.g. \cite{Fr} \cite{HS} \cite{SSu} \cite{Urs}; see \cite{GS3} for a more complete list of references for the untwisted case). 
Here there is also a Chern character map, as part of the data for a differential 
cohomology theory, that lands in a periodic form of de Rham cohomology, where
the periodicity arises from the coefficients of the underlying generalized cohomology theory (see \cite{BNV} \cite{Up}). These differential generalized cohomology theories can in turn be twisted, the foundations of which are given by Bunke and Nikolaus \cite{BN}, building on \cite{BNV} \cite{Urs}.

\medskip
 Carey, Mickelsson, and Wang \cite{CMW} gave a construction 
 of twisted differential K-theory that satisfies the square diagram and short exact sequences. Kahle and Valentino \cite{KV} gave a corresponding list of 
 properties  which can be generalized to any twisted differential   
 cohomology theory. Indeed,  a characterization of twisted differential   cohomology including  twisted differential   K-theory is given in \cite{BN} (except the push-forward axiom).
 The model combines twisted cohomology groups and twisted differential forms in a homotopy theoretic way, analogous to (and as abstract as)  the Hopkins-Singer construction \cite{HS} in the untwisted case.  Having a more concrete geometric description in mind,  Park \cite{Pa} provided a model for twisted differential K-theory in the case that the underlying topological twists represent torsion classes in degree three integral cohomology. 
 Recently, Gorokhovsky and  Lott provided a model of twisted differential   
 K-theory using Hilbert bundles \cite{GL}. 
In this paper, we seek an approach which builds on certain aspects of 
the above works and highlights both homotopic and geometric aspects, 
and which we hope would be  theoretically elegant and conceptually appealing, 
 yet at the same time  being amenable to computations and applications.

\medskip
Among the various spectral sequences one could potentially construct, 
the Atiyah-Hirzebruch spectral sequence (AHSS) would perhaps be the 
most central for twisted theories, because 
of the way the theory is built by gluing data on the underlying space.
As we will explicitly see, this continues to hold for differential twisted 
theories with extra data again assembled from patches on the manifold. 
Hence the {\v C}ech filtration of the manifold by good open covers 
is the most natural in this setting, as it captures the data of 
differential forms, as highlighted classically in \cite{BT}. 
In the topological (and untwisted case) there is an equivalence 
between filtrations on the underlying topological space 
and filtrations on the spectrum of the generalized cohomology theory
\cite{Mau}. While we expect this to hold at the differential 
level \cite{GS3} and also perhaps here at the twisted differential level,
it will become clear that filtrations on the manifold are central to 
the construction and hence are a priori preferred.

\medskip
The AHSS for the simplest nontrivial twisted cohomology theory, namely 
twisted de Rham cohomology, is studied in \cite{AS2}, and more generally
for any odd degree twist in \cite{MW} \cite{LLW}. 
The AHSS for more involved theories has also been studied.
The AHSS in the parametrized setting is described generally but briefly in 
\cite[Remark 20.4.2 and Proposition 22.1.5]{MSi}. 
The twisted K-theory via the AHSS was first discussed in \cite{Ro1} \cite{Ro}
 from the operator algebra point of view, and then briefly in  
 \cite{BCMMS} \cite{AS1} and more extensively in
\cite{AS2} from the topological point of view. 
For $h \in H^3(X; \Z)$, the twisted complex K-theory 
spectral sequence has initial terms with 
$$
E_3^{p,q}=E_2^{p,q}=
\left\{
\begin{array}{ll}
H^p(X; \Z), & q ~{\rm even}, \\
0 & q ~{\rm odd},
\end{array}
\right.
$$
and the first nontrivial differential is the deformation of the cup product 
with twisting class, i.e, 
\footnote{Throughout the paper 
we will use $\tau$ and $\hat{\tau}$ to denote the integral cohomology class 
and its differential refinement, respectively,  for general cohomology theories. 
For K-theory, we will use $h$, $\hat{h}$, and $H$
 to denote, respectively, the 
integral cohomology class, its differential refinement, 
and the corresponding differential form representative.}
$$
d_3=Sq^3_\Z + (-) \cup h: H^p(X; \Z)=E_{3}^{p,q} 
\longrightarrow
 E_{3}^{p+3, q-2}=H^{p+3}(X; \Z)\;.
$$
In the case of twisted cohomology and twisted K-theory the higher 
differentials turn out to be modified by Massey products involving the 
rationalization of the twisting class $h$ \cite{AS2}. For twisted Morava K-theory and twisted Morava E-theory the differentials are discussed in \cite{SW}. All the above theories admit differential refinements and the  AHSS for the corresponding 
untwisted forms of the theories is described in detail by the authors
 in \cite{GS3}. In this paper we construct the AHSS for twisted differential 
 generalized cohomology  theories. Then we focus on  twisted differential K-theory,  whose untwisted  version was considered, 
 among other untwisted theories, in \cite{GS3}.

\medskip
What we are constructing here is the AHSS for twisted differential cohomology
theories which we will denote by $\widehat{\rm AHSS}_{\hat \tau}$, where
$\hat{\tau}$ is a representative of a differential cohomology class, i.e. a higher bundle with connection. When this twisting class if zero in differential cohomology, then
we recover the AHSS for differential cohomology constructed in \cite{GS3}, which
we denote $\widehat{\rm AHSS}$. On the other hand, if we forget the differential 
refinement and strip the theory to its underlying topological theory then 
we recover the AHSS for twisted spectra, which we denote ${\rm AHSS}_{\tau}$, a prominent example of which is that of twisted K-theory 
\cite{Ro1}\cite{Ro}\cite{AS2}. Of course when we take both a trivial twist and no differential 
refinement, then we restrict to the original Atiyah-Hirzebruch case \cite{AH}. 
Summarizing, we will have a correspondence diagram of transformations 
of the corresponding spectral sequences 
$$
\xymatrix{
 \widehat{\rm AHSS}_{\hat \tau} \ar@{~>}[d]_{\hat \tau=0} \ar[rrr]^{I_I} 
 &&& {\rm AHSS}_{\tau} \ar@{~>}[d]^{\tau=0} \\
 \widehat{\rm AHSS} \ar[rrr]^I &&&
 {\rm AHSS}
}
$$
The picture that we have had in mind  for the structure of the corresponding differentials, 
as cohomology operations, is schematically the following 
$$
\xymatrix{
 \widehat{d}_{\hat \tau}=\widehat{\text{Primary}} +
 \widehat{\text{Secondary}} + \cdots
 \ar@{~>}[d]_{\hat \tau=0} \; \ar[rrr]^-{I_I} &&& \;
 {d}_{\tau}= {\text{Primary}} + {\text{Secondary}} + \cdots
  \ar@{~>}[d]^{\tau=0} \\
 \widehat{d}=\widehat{\text{Primary}}
  \ar[rrr]^I &&&
 {d}={\text{Primary}}
}
$$
The upper left corner generalizes the corresponding forms in 
previous works. The identification of $\widehat{d}$ with 
a differential refinement of a primary cohomology operation was
done in \cite{GS2}. The differential secondary refinement 
that will appear in $\widehat{d}_{\hat \tau}$ are essentially 
the differential Massey products constructed in \cite{GS1}. 
The differential $d_\tau$ is that of a twisted cohomology theory, 
as for instance for twisted K-theory, where these secondary operations
are shown to be Massey products (rationally) \cite{AS2}. The bare 
differential $d$ is the Atiyah-Hirzebruch differential which, for the 
case of K-theory, is given by the integral Steenrod square $Sq^3_{\ZZ}$. 
The long term goal that we had in previous projects was to arrive
at the overall picture that we have assembled in this paper. 
In this sense, the current paper can be viewed a culmination of the series  
\cite{GS1}\cite{GS2}\cite{GS3}.

\medskip
Twistings of a suitably multiplicative (that is, $A_\infty$) cohomology theory $R^*$ 
are governed by a space $\op{BGL}_1 (R)$ \cite{MSi}\cite{ABGHR}\cite{ABG}. 
Since twists usually considered arise from cohomology classes, there are maps
 from a source Eilenberg-MacLane space $K(G, m+1)$ for appropriate abelian group 
 $G$ and degree $m$, depending on the theory under consideration. 
We have the following inclusions
$$
\xymatrix{
K(G, m+1) \ar[r] & \op{BGL}_1(R) \; \ar@{^{(}->}[r] & \op{Pic}(R)\;,
}
$$ 
 where $\op{Pic}(R)$ is the $\infty$-groupoid of invertible $R$-module 
 spectra, whose connected component containing the identity 
 is equivalent to $\op{BGL}_1(R)$. For integral twists, 
 i.e. for $G=\Z$, the natural replacement of the Eilenberg-MacLane 
 spaces are  going to be classifying stacks $\mathbf{B}^mU(1)_{\nabla}$.
 We also use the refinement by \cite{BN} of the space 
 $\op{Pic}(R)$. 
 We do so by highlighting the analogies between classical line bundles and $\widehat{K}(X; {\hat{h}})$-theory. We have the following table (see Remark \ref{Rem-table})

\begin{center}
\begin{tabular}{|c||c|}
\hline
 {\bf Line bundles} & {\bf Twisted spectra} 
\\
\hline
\hline
\rule{0pt}{2.5ex}  $\BB\ZZ/2$ & ${\rm Pic}_{\R}^{\rm top}$
\\
\hline 
\rule{0pt}{2.5ex} $\BB\underline{\RR}^{\times}$ & ${\rm Pic}_{\R}^{\rm dR}$
\\
\hline
\rule{0pt}{2.5ex} $\Omega^1_{\rm cl}$ & ${\rm Pic}_{\R}^{\rm form}$
\\
\hline
\rule{0pt}{2.5ex} $\BB\ZZ/2_{\nabla}$ & ${\rm Tw}_{\widehat{\R}}$
\\
\hline
\end{tabular} 
\end{center}
There are classically two viewpoints on real line bundles: 
As bundles with fiber $\RR$ (which is one-dimensional  as an $\RR$-module)
and also as 
locally free invertible ${\cal O}_M$-modules, i.e., 
sheaves which are invertible over the ring of smooth functions.
Similarly, there are two viewpoints on twisted spectra: 
As bundles of spectra ${\cal R}$ (which is rank one as an ${\cal R}$-module),
and also as invertible ${\cal R}$-module spectra. 
The point we highlight is that there is a direct analogy between
each one of the two descriptions of the pair. 
Even more, we can show (see Remark \ref{Rem-equi}) 
that the sheaf cohomology with coefficients in
a line bundle $\cL$ over $M$ is isomorphic to the twisted $H\RR$ cohomology
with twist $\eta: M \to \BB\ZZ/2$ classifying $\cL$.
We will be using stacks, in our context in the sense described in 
\cite{FSSt}\cite{Urs}\cite{SSS3}\cite{cup}. 
In general, twisted topological spectra will satisfy descent over
the base space, hence are stacks \cite{BN}.  For the case of
 twisted topological K-theory this is also shown in \cite{AW}.
Differential theories similarly satisfy descent over 
the base manifold $M$ \cite{BN}.


\medskip
The paper is organized as follows.
Section \ref{Sec-gen} we provide our slightly modified take on \cite{BN} 
on Twisted differential spectra for the purpose of constructing 
 their AHSS.  
 In particular, in Section  \ref{Sec-rev} we start recalling the construction of
 \cite{BN} with slight modifications to suit our purposes. 
 This leads us in Section  \ref{Sec-can} to 
  describe a canonical bundle of spectra over the Picard $\infty$-groupoid 
  which comes from the $\infty$-Grothendiek construction. 
 In Theorem \ref{Thm-equiv} we present several different, but equivalent, ways to think about twisted differential cohomology theories. 
 We use this in  Section  \ref{Sec-flat}  to provide 
 analogies between these objects and the notion of a line bundle equipped with flat connection, leading to the identification of  {\it a twisted differential theory as a bundle of spectra equipped with flat connection}, whose meaning we explain. 
 This is captured
 in the following table, continuing the analogies in the above table. 
 
 \vspace{5mm}
\begin{table}[!hp]
\label{table}
\hspace{-.5cm}
\footnotesize{
\begin{tabular}{|c||c|}
 \hline
  \rule{0pt}{4ex} 
{\large \bf Line bundle }
& 
{\large \bf Bundle of spectra}
\rule[-2ex]{0pt}{0pt} 
\\
\hline
\hline
 \rule{0pt}{4ex} 
$$
\xymatrix{\RR  \ar[r]\ar[d]& E\ar[d]^-{\pi}
\\
\ast ~\ar@{^{(}->}[r] & M 
}
$$
& 

$$
\xymatrix{H\big(  \RR[u, u^{-1}]\big) 
\cong_{\rm ch} K \wedge H\RR \ar[d] \ar[r]& E\ar[d]^-{\pi}
\\
\ast ~\ar@{^{(}->}[r]  & M 
}
$$
\\
\hline
 \rule{0pt}{4ex} 
{\bf Set of global sections }
&
{\bf Space of global sections}
\\
$
\xymatrix{
\Gamma(M, E)\cong  \lim \{  \prod C^\infty (U_\a, \RR)
\ar@<.05cm>[r]^-{g_{\a \b}} 
 \ar@<-.05cm>[r] & \prod C^\infty (U_{\a \b}, \RR)
\}
}
$
&
$
\Gamma(M, E)\cong  \xymatrix{ \op{holim} \{  \prod \Omega^* (U_\a)
\ar@<.05cm>[r]^-{e^{dA_{\a \b}}}
 \ar@<-.05cm>[r] & \prod \Omega^* (U_{\a \b})
\ar@<.1cm>[r]^-{CS}\ar[r] \ar@<-.1cm>[r] &
\cdots
\}
}
$
\rule[-4ex]{0pt}{0pt} 
\\
\hline
 \rule{0pt}{4ex} 
{\bf Local ring of smooth real functions}
&
{\bf Local DGA of smooth functions}
\\

$
\Omega^* \otimes \RR[0] \simeq C^\infty(-, \RR)
$
&
$
\Omega^* \otimes \RR[u, u^{-1}]\simeq \{ 
\Omega^{\rm ev} \longrightarrow \Omega^{\rm odd} \longrightarrow 
\cdots \}
$
\\
&

\\
\hline
 \rule{0pt}{4ex} 
{\bf Flat connection}
&
{\bf Flat connection}
\\
$
\nabla: \Gamma(M, E) \longrightarrow \Gamma(M, E \otimes T^*M)
$
&
$
d_H: \Gamma(M, E) \longrightarrow \Gamma(M, E \wedge \Exterior^3(T^*M))
$
\\
$\nabla(\sigma \cdot s)= d\sigma s + \sigma \nabla s$
\; with \; $\nabla^2=0$
&

$d_H(\omega \wedge \alpha)=d_H(\omega) \wedge \alpha + \omega \wedge d\a$
\; with \; $d_H^2=0$
\rule[-3ex]{0pt}{0pt} 
\\
\hline
 \rule{0pt}{4ex} 
{\bf Parallel section}
&
{\bf Parallel section}
\\
$s:M\to E$, \; $ \nabla s=0$
&
$d_H$-parallel: $d_H(\omega)=0$ 
\rule[-3ex]{0pt}{0pt} 
\\
\hline
 \rule{0pt}{4ex} 
{\bf Discrete transition functions}
&
{\bf Topological twist}
\\ 
$g_{\a \b}\in \check{H}^1(M;\RR^{\times}) $
& 
$\eta_{\a\b\gamma\delta}\in \check{H}^3(M;\ZZ)$
\\
written in local trivializations: 
& 
written in local trivializations: 
\\
$s_\a:U_{\alpha}\to E$,  \;$\nabla(s_\a)=0$
&
$e^{B_{\alpha}}\wedge \omega_{\alpha}$, \; 
$d_{H}(e^{B_{\alpha}}\wedge \omega_{\alpha})=0$
\rule[-3ex]{0pt}{0pt} 
\\
\hline
 \rule{0pt}{4ex} 
{\bf  Representation}
 &
 {\bf Representation}
 \\
Action of $\RR^{\times}$ on $\RR$ as multiplication by units
 &
Action of $BU(1)$ on $K$ via line bundles.
\rule[-3ex]{0pt}{0pt} 
\\
\hline
\end{tabular} 
}
\end{table}

\vspace{.5cm}
 The general construction of $\widehat{\rm AHSS}_{\hat \tau}$,
 the AHSS for any twisted differential cohomology theory, is 
given in Section  \ref{Sec-AHSS}. The differentials are then 
identified on general grounds and the corresponding properties
are discussed in Section  \ref{Sec-prop}.
Having provided a general discussion, we start focusing on 
twisted differential K-theory in Section  \ref{Sec-K}. 
Indeed, in Section  \ref{Sec-gerb}
we give a detailed description of the twists in differential $K$-theory 
from various perspectives. 
In order to arrive at explicit formulas, we describe the 
Chern character in twisted differential K-theory in Section  \ref{Sec-ch}.
Then in the following two sections we characterize the 
differentials more explicitly. Unlike the classical AHSS, in the 
differential setting, as  also noticed for the untwisted case \cite{GS3}, we find
 two kinds of differentials:  The low degree differentials are
characterized in  Section  \ref{Sec-low}, while the higher flat ones  
are described in Section  \ref{Sec-high}. We end with explicit and detailed 
examples, illustrating the machinery, in Section  \ref{Sec-ex}.
We calculate twisted differential K-theory of the 3-sphere
(or the group $\op{SU}(2)$)
$\widehat{K}^*(S^3; {\hat{h}})$ in two ways: Using the 
Mayer-Vietoris sequence  of Proposition \ref{prop twdiff}(iii)
and using the $\widehat{\rm AHSS}_{\hat \tau}$ of Theorem 
\ref{Thm-ourAHSS}, as adapted to K-theory with differential twist 
$\hat{h}$ in previous sections. 
While we present both approaches as being useful,  
the highlight is that the latter approach is much more powerful,
at least in this case.

 \newpage

%
%
%

\section{Twisted differential spectra and their AHSS} 
\label{Sec-gen}

\subsection{Review of the general construction of differential twisted spectra}
\label{Sec-rev}

In this section we start by recounting some of the main constructions 
in \cite{BN} and then provide a description of local triviality as
well as general properties of twisted differential spectra. 
The point of view taken in \cite{BN} uses primarily Picard $\infty$-groupoids to define the differential twists of a differential function spectrum. Topologically, given an $E_{\infty}$ (or even an $A_\infty$) ring spectrum $\R$, one can obtain this $\infty$-groupoid as a nonconnected delooping of its $E_{\infty}$ space of units. More precisely, if ${\rm Pic}_{\R}$ denotes the infinity groupoid of invertible $\R$-module spectra (with respect to the smash product $\wedge_{R}$), then we have an equivalence
$$
\op{B GL}_1(\R)\times \pi_0\big({\rm Pic}_{\R}\big) \simeq {\rm Pic}_{\R}\;,
$$
where $\pi_0\big({\rm Pic}_{\R}\big)$ is the group of isomorphism classes of invertible $\R$-module spectra (see \cite{ABGHR} \cite{ABG} \cite{SW}).  

\medskip
In order to discuss the differential refinement of twisted cohomology, we need to retain some geometric data which cannot be encoded in spaces. Indeed, twisted differential cohomology theories crucially use a sort of twisted de Rham theorem to connect information about differential forms on a smooth manifold with some cohomology theory evaluated on the manifold. This mixture of geometry and topology is effectively captured by smooth stacks, with the homotopy theory being captured by the simplicial direction of the stack and the geometry being captured by the local information encoded in the site. In general, we would like to consider not only $\infty$-groupoid valued smooth stacks, but also stacks valued in the stable $\infty$-category of spectra; see \cite{Lur}\cite{Lur2} for more comprehensive accounts of these constructions. 
\footnote{
The reader who is not interested in the technical details can regard this rigorous development as 
 a sort of black box, bearing in mind that in the language of $\infty$-categories no commutative  diagram  is strict: only commutative up to higher homotopy coherence. The familiar concepts  at the 1-categorical level generally hold in the $\infty$-context, but only up to homotopy coherence.
 }

\medskip
Purely topological theories can be regarded as \emph{constant} smooth theories. For example, let $\R$ be an `ordinary' ring spectrum. Then we can consider the assignment which associates to each manifold $M$ the spectrum $\R$. This assignment defines a prestack, i.e., a functor
$$
{\rm const}(\R):\mf \longrightarrow \Sp
$$ 
on the category of smooth manifolds with values in spectra. 
We can equip the category $\mf$ with the Grothendieck topology of good open covers, turning it into a site. The stackification of ${\rm const}(\R)$ with respect to this topology will be denoted $\underline{\R}$. It is obviously locally constant,
i.e. constant on some open neighborhood of each point of $M$. 

\medskip
 For a fixed manifold $M$, we can restrict the above functor ${\rm const}(\R)$ to the overcategory  $\mf/M$ (i.e. the $\infty$-category on arrows $N\to M$), equipped with the restricted coverage.  After stackification, this gives a locally constant sheaf of spectra $\underline{\R}$ on $M$. More generally, we can do this for any locally constant sheaf of spectra on $\mf$. Then we can consider the monoidal $\infty$-category of \emph{locally constant} sheaves of $\underline{\R}$-module spectra over $M$ and consider the space of invertible objects. \footnote{Here the monoidal structure is given by the $R$-module smash product $\wedge_{R}$. Clearly we are omitting the details of monoidal  $\infty$-categories and the symmetric monoidal smash product. It would be far too lengthy to  include these details here and we encourage the reader to consult \cite{Lur2} for more information. }  
  In \cite{BN}, this $\infty$-groupoid is denoted by ${\rm Pic}_{\underline{\R}}(M)$. In order to simplify notation and to highlight the fact that this object is related to the twists of the underling theory $\R$, we will denote this $\infty$-groupoid by ${\rm Pic}_{\R}^{\rm top}(M)$. This is further motivated by the following properties \cite[Sec. 3]{BN}:

\begin{enumerate}\label{prop loc sp}

\vspace{-2mm}
\item {\bf (Descent):} The assignment $M\mapsto {\rm Pic}^{\rm top}_{\R}(M)$ defines a smooth stack on the site of all manifolds, topologized with good open covers.

\vspace{-2mm}
\item {\bf (Correspondence with topological twist):}
For a manifold $M$, we have equivalences ${\rm Pic}^{\rm top}_{\R}(M)\simeq \map(M,{\rm Pic}_{\R})\simeq \underline{{\rm Pic}_{\R}}(M)$ with ${\rm Pic}_{\R}$ the usual Picard space of invertible $\R$-module spectra, the second mapping $\infty$-groupoid is the nerve of the mapping space in CGWH-spaces and $\underline{{\rm Pic}_{\R}}$ is the constant stack associated to the nerve of ${\rm Pic}_{\R}$. All three of the resulting spaces model the $\infty$-groupoid of invertible $\R$-modules parametrized over $M$.

\vspace{-2mm}
\item {\bf (Cohomology theory from spectrum):} Let $M$ be a manifold with $\R$-twist ${\R}_{\tau} \in \Pi_1{\rm Pic}_{\R}^{\rm top}(M)$, where $\Pi_1$ denotes the fundamental groupoid. Then the classical $\E$-twisted cohomology of $M$ is computed via \footnote{Generally, we will denote a spectrum by a script font (e.g. $\R$), while the underlying theory it represents will be denoted by the standard math display font (e.g. $R$).}
$$
R(M ; {\tau}):=\pi_0(\Gamma(M, {\R}_{\tau}))\;,
$$
where $\Gamma$ is the global sections functor $\Gamma:\sh_{\infty}(M;\Sp)\to \Sp$.
\end{enumerate}

\begin{remark}
[Determinantal twist vs. dimension shift]
\label{Rem-dim}
Notice that the third property does not include a degree of the theory. This is because shifts of the spectrum $\tau$ are also $\R$-twists. Thus the various degrees of 
$R(-;{\tau})$ can be obtained by considering the shifted $\R$-twists. However, because we are  dealing with differential cohomology
and are interested in geometry and applications, we find it more practical to 
separate the two notions. Thus, if we have a combination of twists of the form $\tau$ and ${\bf n}$, where ${\bf n}$ is the twist given by the dimension 
shift $\Sigma^n\R$ we will write
$$
R^n(M; {\tau}):= \pi_0 (\Gamma (M; \R_{\tau+{\bf n}}))
$$
for the cohomology of the product 
$\R_{\tau}\wedge_{\R}\Sigma^n\R\simeq \Sigma^n\R_{\tau}$.

\end{remark}

The above discussion does not yet taken into account differential data. It is merely the manifestation of the underlying topological theory, embedded as a geometrically discrete object in smooth stacks. To describe the differential theory, we first need to discuss the differential refinement of an $E_{\infty}$-ring spectrum $\R$. This is discussed extensively in \cite{BN} and we review only the essential pieces discussed there, adapted to our purposes. 

\begin{definition}
Define the \emph{category of differential ring spectra}
 as the $(\infty,1)$-pullback
$$
\xymatrix@=1.5em{
{\rm diff}(\mathscr{R}{\rm ing}\mathscr{S}{\rm p})\ar[d] \ar[rr] && 
\RR\hbox{-} {\rm CDGA}\ar[d]^{H}
\\
\mathscr{R}{\rm ing}\mathscr{S}{\rm p} \ar[rr]^-{\wedge H\RR} && 
\RR\hbox{-}\mathscr{A}{\rm lg}\mathscr{S}{\rm p}\;.
}
$$
Here $\mathscr{R}{\rm ing}\mathscr{S}{\rm p}$ denotes the $(\infty,1)$-category of $E_{\infty}$-ring spectra, 
$\RR$-$\mathscr{A}{\rm lg}\mathscr{S}{\rm p}$ is the $(\infty,1)$-category of commutative $\RR$-algebra spectra, and $\RR$-${\rm CDGA}$ is the $(\infty,1)$-category of commutative differential graded algebras. The functor $H$ is the Eilenberg-MacLane functor and $\wedge H\RR$ takes the smash product with the real Eilenberg-MacLane spectrum (see, e.g.,  \cite{Sh} for details).
\end{definition}

\begin{remark}
\noindent {\bf (i)} Note that the objects of the $(\infty,1)$-category ${\rm diff}(\mathscr{R}{\rm ing}\mathscr{S}{\rm p})$ can be identified with triples $(\R,c,A)$, where $\R$ is an $E_{\infty}$-ring spectrum, $A$ is a ${\rm CDGA}$ and $c$ is an equivalence of ring spectra
$$
c:\R \wedge H \RR\overset{\simeq}{\longrightarrow} HA\;.
$$
\item {\bf (ii)} We also note that it was shown in \cite{Sh} that the Eilenberg-MacLane functor is essentially surjective onto $H\ZZ$-module spectra, from which we can deduce that every ring spectrum admits a differential refinement since $H\RR$-module spectra can be regarded as $H\ZZ$-module spectra. 
\end{remark}

\begin{definition}
\label{def-dRA}
 Let $A$ be a ${\rm CDGA}$. We define the \emph{de Rham complex with values in $A$} to be the sheaf of chain complexes
 $$
 \Omega^*(-;A):=\Omega^*\otimes_{\RR} \underline{A}\;,
 $$
 where $\underline{A}$ is the locally constant sheaf of ${\rm CDGA}$'s associated to $A$.
 \end{definition}
 To define the smooth stack of twists, Bunke and Nikolaus \cite{BN} 
 define a certain stack which contains all the twisted de Rham complexes 
 corresponding to the theory. We now recall this stack.
Let ${\cal L}$ be a sheaf of $\Omega^*(-;A)$-modules on $M$ (i.e. on the restricted site $\mf/M$). Then
\begin{enumerate}

\vspace{-2mm}
\item ${\cal L}$ is  called \emph{invertible} if there is ${\cal K}$ such that ${\cal L}\otimes_{\Omega^*(-;A)}{\cal K}$ is isomorphic to $\Omega^*(-;A)$.

\vspace{-3mm}
\item  ${\cal L}$ is called $K$-\emph{flat} if the functor
$$
{\cal L}\otimes_{\Omega^*(-;A)}-:\Omega^*(-;A)\hbox{-}\mathscr{M}{\rm od}
\longrightarrow
 \Omega^*(-;A)\hbox{-}\mathscr{M}{\rm od}
 $$
preserves object-wise quasi-isomorphisms.

\vspace{-3mm}
\item ${\cal L}$ is called \emph{weakly locally constant}, if it is quasi-isomorphic to a locally constant sheaf of modules.
\end{enumerate}

\noindent 
Denote by $\op{Pic}^{\rm form}$
the stack which assigns to each smooth manifold $M$, the Picard groupoid 
$\op{Pic}^{\rm form}(M)$
of all invertible, K-flat, weakly locally constant sheaves of $\Omega^*(-;A)$-modules. The passage from the topological twists to differential twists is then obtained via the $(\infty,1)$-pullback in smooth stacks
\(
\label{bun twpull}
\xymatrix@=1.6em{
{\rm Tw}_{\widehat{\R}} \ar[rrr]\ar[d] &&&
 {\rm Pic}^{\rm form}\ar[d]^-{H}
\\
{\rm Pic}^{\rm top} \ar[rrr]^-{\wedge H\RR} &&&
 {\rm Pic}^{\rm dR}\;,
}
\)
where $\op{Pic}^{\rm top}$, $\op{Pic}^{\rm form}$, and $\op{Pic}^{\rm dR}$ have the 
more detailed respective notations ${\rm Pic}^{\rm loc}_{\underline{\R}}$,
${\rm Pic}_{\Omega^*(-;A)}^{\rm wloc,fl}$, and  ${\rm Pic}^{\rm loc}_{H(\Omega^*(-;A))}$
in \cite{BN}.

\medskip
A bit of explanation is in order. First, the condition of weakly locally constant for the 
sheaf of $\Omega^*(-;A)$-modules ${\cal L}$ gives rise to a twisted de Rham theorem. Indeed, for each fixed point $x\in M$, there is an open set $x\in U\subset M$ such that ${\cal L}(U)$ is quasi-isomorphic to a constant complex on $U$. This locally constant complex represents a cohomology theory with local coefficients on $M$ and these local quasi-isomorphisms induce a quasi-isomorphism of sheaves of complexes on $M$. At the level of cohomology, the induced map can be regarded as a twisted de Rham isomorphism.  Next, the condition of $K$-flatness ensures that the tensor product is equivalent to the derived tensor product and therefore preserves invertible objects. This way, we have a well-defined Picard groupoid. 

\begin{remark}
Note that an object in the $(\infty,1)$-category ${\rm Tw}_{\widehat{\R}}$ is given by a triple $({\R}_{\tau}, t, {\cal L})$, with ${\R}_{\tau}$ a topological $\R$-twist, ${\cal L}$ an invertible, K-flat, weakly locally constant module over $\Omega^*(-;A)$ and $t$ a zig-zag of weak equivalences (or a single equivalence in the localization at quasi-isomorphisms)
$$
t:H(\L)\overset{\simeq}{\longrightarrow} {\R}_{\tau} \wedge H\RR\;.
$$
The equivalence $t$ exhibits a twisted de Rham theorem for the rationalization of the topological twist ${\R}_{\tau}$. 
\end{remark}

\begin{definition}
Such a triple $({\R}_{\tau}, t, {\cal L})$ is called a \emph{differential refinement  of a topological twist} ${\R}_{\tau}$.
\end{definition}

One of the hallmark features of twisted cohomology is that of local triviality. Indeed, if we are to think of twisted cohomology theories as a bundle of theories, parametrized over some base space (in the sense of \cite{MSi}), then one expects that if we restrict to a contractible open set on the parametrizing space, the theories should all be isomorphic to the `typical fiber', given by the underlying theory which is being twisted. This happens for differential cohomology theories, just as it does for topological theories. 

\begin{proposition}
[Local triviality of the twist]
\label{prop-triviali}
Let $\widehat{\R}_{\hat{\tau}}:=(\R_{\tau},t,{\cal L})$ be a sheaf of spectra refining the twisted theory $\R_{\tau}$ on a smooth manifold $M$. Suppose the underlying topological twist comes from a map 
$
\tau: M \to  B{\rm GL}_1(\R) \into  \op{Pic}_{\R}^{\rm top}
$.
Then, around each point $x\in M$, there is an open set $U$ such that, when restricted 
via the change of base functor
$$
i^*:{\rm Tw}_{\widehat{\R}}(M)\longrightarrow {\rm Tw}_{\widehat{\R}}(U)\;,
$$
induced by the inclusion $i:U\to M$, we have an equivalence in ${\rm Tw}_{\widehat{\R}}(U)$
$$
\widehat{\R}_{\hat{\tau}}\simeq \widehat{\R}\;,
$$
where $\widehat{\R}=(\R,c,A)$ is a differential refinement of the untwisted theory $\R$.
\end{proposition}
\theproof
Let $x:\ast\to U\subset M$, be a point in $M$. 
Since $M$ is a manifold, without loss of generality, it suffices to prove the claim with $U$ contractible.  In this case, the restriction of $\R_{\tau}$ to $U$ trivializes and we have an equivalence $\R_{\tau}\simeq \R$
in $\sh_{\infty}(U;\Sp)$. As part of the data for the differential refinement of $\R$, we have an equivalence in $\sh_{\infty}(\mf, \Sp)$
$$
H(\Omega^*(-;A))\simeq \R\wedge H\RR\;,
$$
and hence an equivalence on the restricted site $\sh_{\infty}(U;\Sp)$. Combining, we have equivalences in $\sh_{\infty}(U;\Sp)$,
$$
\xymatrix@R=1.8em{
\R \ar[rr]^\simeq \ar[d]_-{\wedge H \RR} &&
\R_\tau \ar[d]^-{\wedge H \RR}
\\
\R\wedge H\RR \ar[d]_{\simeq}^{c} \ar[rr]^\simeq &&
\R_\tau\wedge H\RR \ar[d]^{\simeq}_{t} 
\\
H(\Omega^*(-;A)) \ar[rr]^{\simeq } && H({\cal L})\;,
}
$$
where the bottom map in the diagram is defined by composing 
the evident three equivalences.  
By the definition of localization and the properties of $H$, this implies that $\Omega^*(-;A)$ and ${\cal L}$ are quasi-isomorphic. Putting this data together, we see that we have defined an equivalence $\widehat{\R}_{\hat{\tau}}\overset{\simeq}{\to} \widehat{\R}$ in the $\infty$-groupoid ${\rm Tw}_{\widehat{\R}}(U)$.
\endofproof

We now summarize some of the properties of twisted differential cohomology
from above.  

\begin{proposition}
[Properties of twisted differential cohomology]
\label{prop twdiff}
Let $\widehat{R}^*(-; {\hat{\tau}})$ be a twisted counterpart of a differential cohomology theory $\widehat{R}^*$ on a smooth manifold $M$. Then 
$\widehat{R}^*(-; {\hat{\tau}})$ satisfies the following:
\begin{enumerate}[label=(\roman*)]

\vspace{-2mm}
\item \emph{(Additivity)} For a disjoint union $f:N=\coprod_{\alpha} N_{\alpha}\to M$ of maps $f_{\alpha}:N_{\alpha}\to M$, we have a decomposition
$$
\widehat{R}^*(N; {\hat{\tau}})\cong \bigoplus_{\alpha}
\widehat{R}^*(N_{\alpha}; {\hat{\tau}})\;.
$$
\vspace{-6mm}
 \item \emph{(Local triviality)} If $\hat{\tau}$ has underlying topological twist $\tau$ coming from a map $\tau:M\to BGL_1(\R)$, then for every point $x\in M$, there is an open neighborhood $x\in U\overset{i}{\into}M$ such that the restriction of $\widehat{R}^*(-; {\hat{\tau}})$ to $U$ trivializes. That is we have a natural isomorphism of functors 
$$
i^* \widehat{R}^*(-; {\hat{\tau}})\cong \widehat{R}(-)\;,
$$
on the over category $\mf/U$.
\vspace{-2mm}
\item \emph{(Mayer-Vietoris)} Let $\{U,V\}$ be an open cover of $M$. Then 
we have a long exact sequence
\begin{equation*}
\xymatrix{
  &
 ... \ar[r] &  
{\R^{-2}_{\RR/\ZZ}(U; \hat{\tau})
\oplus
  \R^{-2}_{\RR/\ZZ}(V; \hat{\tau})}
 \ar[r]
 & 
 \R^{-2}_{\RR/\ZZ}({U\cap V}; \hat{\tau})
 \ar`r[d]`[l]`[llld]`[dll][dll]  \\
  &\widehat \R^{0}(M; \hat{\tau}) \ar[r] & 
  \widehat \R^0(U; \hat{\tau})\oplus \widehat R^0(V; \hat{\tau}) \ar[r]& 
\widehat \R^0(U\cap V;  \hat{\tau}) \ar`r[d]`[l]`[llld]`[dll][dll] 
\\ 
 & \R^{1}(M; {\tau}) \ar[r] &   \R^{1}(U; {\tau}) 
 \oplus \R^{1}(V; {\tau})  \ar[r]& ...
 }
\end{equation*}
\end{enumerate}
\end{proposition}

\theproof
The first follows from the fact that the global sections on a coproduct is isomorphic to the product of global sections and that taking connected components commutes with products. The second property arises from Proposition \ref{prop-triviali}.
 The third property is stated in \cite[Prop. 5.2]{BN}, but we give 
 a detailed proof here. With respect to the coverage on $\mf/M$ via good open covers, the sheaf of spectra $\widehat{\R}_{\hat{\tau}}$ satisfies descent. Consider the pushout diagram
 $$
 \xymatrix{
 U\cap V\ar@<.1cm>[r]\ar@<-.1cm>[r] & U\coprod V \ar[r] & U\cup V =M
 }\;.
 $$
Since $\widehat{\R}_{\hat{\tau}}$ satisfies descent on $M$, the pullbacks to $U\coprod V $ and $U\cap V$ also satisfy descent and (by restriction) yield an exact sequence 
$$
\widehat{\R}_{\hat{\tau}}\vert_{U\cap V}\to \widehat{\R}_{\hat{\tau}}\vert_{U}\times \widehat{\R}_{\hat{\tau}}\vert_{V}\to \widehat{\R}_{\hat{\tau}}\vert_{U\cap V}
$$
of sheaves of spectra. After taking global sections, we get the iterated fiber/cofiber sequence
 $$
 \hdots \to \Sigma^{-1}\widehat{\R}_{\hat{\tau}}(M)\to \widehat{\R}_{\hat{\tau}}(U\cap V)\to \widehat{\R}_{\hat{\tau}}(U)\times \widehat{\R}_{\hat{\tau}}(V) \to  \widehat{\R}_{\hat{\tau}}(M)\to \hdots \;.
 $$
 Finally, passing to connected components, one identifies the long exact sequence as stated.
\endofproof

\medskip
The local triviality property will be useful throughout and  
 the Mayer-Vietoris property will be explicitly utilized
  in the examples in  Sec. \ref{Sec-ex}. 

\medskip
Now we have seen that a map $\hat{\tau}:M\to {\rm Tw}_{\hat{\R}}$ specifies a twist $\widehat{\R}_{\hat{\tau}}$ for the differential ring spectrum $\widehat{\R}$ and we know how to define the twisted cohomology abstractly via global sections 
of this sheaf of spectra. In practice, however, it is useful to understand how the 
sheaf $\widehat{\R}_{\tau}$ is built out of local data, just as it is useful to understand the transition functions of a vector bundle. Happily, descent allows us to understand how a map to the stack of twists specifies local gluing data via pullback by a universal bundle over the twists, although we next need to develop a bit of machinery to make this precise.

\subsection{The canonical bundle associated to the Picard $\infty$-groupoid}
\label{Sec-can}

In this subsection, we describe a canonical bundle of spectra over the Picard $\infty$-groupoid which comes from the $\infty$-Grothendiek construction. While rather abstract at first glance, we will see that this perspective has many conceptual and computational advantages. The general framework for this construction was set up by Jacob Lurie in \cite{Lur2}. Here we are simply unpackaging and 
adapting the general construction to our context.
 
\medskip
In \cite[Theorem 3.2.0.1]{Lur2}, a general Grothendiek construction for $\infty$-categories is described for each simplicial set $S$, whose associated simplicial category is equivalent to a fixed simplicial category $\C^{\rm op}$. The construction gives a Quillen equivalence (depending on a choice of equivalence $\phi:S\to \C^{\rm op}$) 
$$
\xymatrix{
\op{St}^{+}_{\phi}:  ({\sset}^+)_{/S}
\; \; \ar@<.1cm>[r] & \ar@<.1cm>[l] \; \;
(\sset^+)^{\C}: {\op{Un}^+_{\phi}}
}
\;,$$
where the right hand side is a presentation of the $\infty$-category of functors ${\rm Fun}(S,\mathscr{C}{\rm at}_{\infty})$. 
\footnote{The functors $\op{St}^{+}_{\phi}$
and ${\op{Un}^+_{\phi}}$ are called the straightening and unstraightening 
functors, respectively.  However, we will not need anything explicit 
about such an interpretation here.}
The fibrant objects on the left hand side are precisely Cartesian fibrations $X\to S$. Thus the construction associates to each functor 
$\overline{p}:\C\to \mathscr{C}{\rm at}_{\infty}$ 
\footnote{where $\mathscr{C}{\rm at}_{\infty}$ is the $\infty$-category of 
all $\infty$-categories.}
a Cartesian fibration $p:X\to S$ whose fiber at $\phi^{-1}(c)$ can be identified with the $\infty$-category $\overline{p}(c)$. 
\footnote{The analogy to keep in mind here is the tautological  line 
bundle over $\CC P^n$, each point  of which is represented by
a line in $\CC^{n+1}$, and the associated bundle  is then defined
by assigning to each one of these points a copy of that line as the 
fiber.}

\medskip
Now, in all 
$\infty$-categories, i.e., in $\mathscr{C}{\rm at}_{\infty}$
there are the stable ones
$$
\xymatrix{
i:\mathscr{S}{\rm tab}_{\infty}\; \ar@{^{(}->}[r] &
 \mathscr{C}{\rm at}_{\infty}
}\;,
$$
and this inclusion functor admits a right adjoint, ${\rm Stab}$, which stabilizes an $\infty$-category. In particular, we can apply this functor to the codomain fibration $\C^{\Delta[1]}\to \C$, which gives the \emph{tangent} $\infty$-category $p:T{\C}\to \C$. This map is a Cartesian fibration, and the associated functor under the Grothendiek construction sends each $c\in \C^{\rm op}$ to the stabilization of the slice ${\rm Stab}(\C_{/c})$ (see \cite[Sec. 1.1]{Def}). In particular, if we take $\C=\infty\mathscr{G}{\rm rpd}$ and fix an $\infty$-groupoid $X\in \infty\mathscr{G}{\rm rpd}$, again by the $\infty$-Grothendiek construction, we get an equivalence
$$
 {\rm Stab}(\infty\mathscr{G}{\rm rpd}_{/X}) \simeq 
 {\rm Stab}({\rm Fun}(X,\infty\mathscr{G}{\rm rpd}))\;.
$$
By the discussion above, the stabilization on the left is precisely the fiber of the tangent $\infty$-category $T({\infty\mathscr{G}{\rm pd}})\to \infty\mathscr{G}{\rm pd}$ at $X$. From the stable Giraud theorem (see \cite[Remark 1.2]{Lur}), one sees that the right hand side is equivalent to ${\rm Fun}(X,\mathscr{S}{\rm p})$. 

\medskip
Now bringing all this down to earth, what we have is that a functor $\overline{p}:X\to \mathscr{S}{\rm p}$ is equivalently an object in the tangent $\infty$-category at $X$. Such an object can be expressed as a `map' 
\footnote{The main point here is that $E$ and $X$ live in different categories. 
However, we can get an actual map if were to embed $X$ into its 
tangent $\infty$-category as the terminal object. 
For ordinary $\infty$-groupoids, this sends $X$ to 
the corresponding spectrum  $\Sigma^\infty X_+$.}
$p:E\to X$, which is to be thought of as a bundle of spectra over $X$. Notice also that one can do usual bundle constructions in the tangent $\infty$-category. In particular, given a map $f:X\to Y$ of $\infty$-groupoids, the change-of-base functor assigns to each object $E\to Y$ in $T({\infty\mathscr{G}{\rm pd}})_{Y}$ 
 an object $f^*(E)\to X$ in $T({\infty\mathscr{G}{\rm pd}})_{X}$. This construction gives a pullback diagram in the global tangent $\infty$-category $T({\infty\mathscr{G}{\rm rpd}})$ given by
$$
\xymatrix@R=1.5em{
f^*(E) \ar[d] \ar[rr] && E\ar[d] 
\\
X \ar[rr]^-{f} && Y
\;,}
$$
where $X$ and $Y$ are embedded via the ``global sections functor" $s:\infty\mathscr{G}{\rm rpd}\to T({\infty\mathscr{G}{\rm rpd}})$, which assigns each infinity groupoid $X$ to the constant functor $X:\Delta[1]\to \infty\mathscr{G}{\rm rpd}$ and then to the corresponding object in the stabilization. 

\medskip
In our case, we take $X$ to be the Picard $\infty$-groupoid ${\rm Pic}_{\R}^{\rm top}$ (see the description just before Remark \ref{Rem-dim}) and define the functor 
$$
\overline{p}:{\rm Pic}_{\R}^{\rm top}\longrightarrow
 \mathscr{S}{\rm p}
$$
as the functor which assigns to each module spectrum $\R_{\tau}$ over $\R$
the spectrum $\R_{\tau}$. For tautological reasons this assignment defines an $\infty$-functor and we denote the associated bundle of spectra by $\pi:\lambda^{\rm top}_{\R}\to {\rm Pic}_{\R}^{\rm top}$. Since a map $\tau:\ast\to {\rm Pic}^{\rm top}_{\R}$ is nothing but a choice of object of ${\rm Pic}^{\rm top}_{\R}$, the construction immediately identifies the fiber with $\R_{\tau}$. That is, we have 
the following diagram in the tangent $\infty$-category
$$
\xymatrix@R=1.5em{
\R_{\tau} \ar[d] \ar[rr] && \ar[d] \lambda^{\rm top}_{\R}
\\
\ast \ar[rr]^-{\tau} && {\rm Pic}^{\rm top}_{\R} \;.
}
$$
One pleasant feature about this point of view is that all the machinery holds equally well for any $\infty$-topos. In particular, if we replace $\infty\mathscr{G}{\rm pd}$ everywhere by $\sh_{\infty}(M)$, the $\infty$-topos of smooth stacks on $M$,  then all the constructions hold equally well. We can, therefore, consider the functor 
\(
\label{pm func}
\overline{p}_{M}:{\rm Tw}_{\widehat{\R}}(M)\longrightarrow
 \sh_{\infty}(M;\mathscr{S}{\rm p})\;,
\)
which assigns to each twist $\widehat{\R}_{\hat{\tau}}$ the sheaf of spectra $\widehat{\R}_{\hat{\tau}}$ on $M$. Associated to this functor, there is an associated bundle of (now) sheaves of spectra $p_{M}:\lambda_{\widehat{\R}}(M)\to {\rm Tw}_{\widehat{\R}}(M)$ and the fiber at a point $\tau:\ast\to {\rm Tw}_{\widehat{\R}}(M)$ is the twist $\widehat{\R}_{\hat{\tau}}$ on $M$.
\begin{definition}
[Canonical bundle of sheaves of spectra]
We define the \emph{canonical bundle of sheaves of spectra} over the stack of twists ${\rm Tw}_{\widehat{\R}}$ as the presheaf with values in $T(\infty\mathscr{G}{\rm rpd})$ which associates to each smooth manifold $M$ the object 
$$
p_{M}:\lambda_{\widehat{\R}}(M)\longrightarrow {\rm Tw}_{\widehat{\R}}(M)
$$
associated to the tautological functor 
$\overline{p}_{M}$ given in \eqref{pm func} 
which assigns to each object $\widehat{\R}_{\hat{\tau}}$ the sheaf of spectra $\widehat{\R}_{\hat{\tau}}$ on $M$. We denote this object by $\lambda_{\widehat{\R}}\to {\rm Tw}_{\widehat{\R}}$.
\end{definition} 

Essentially, because the stabilization of an $\infty$-category is built as a limit of $\infty$-categories, we have a canonical equivalence
$$
\mathscr{PS}{\rm h}_{\infty}(\mf,T(\infty\mathscr{G}{\rm rpd}))\simeq T(\mathscr{PS}{\rm h}_{\infty}(\mf))\;.
$$
By manipulating adjoints and using the fact that the stabilization functor preserves limits, we see that this equivalence holds even at the level of sheaves. This implies that (under the Grothendiek construction) the bundle $p_{M}$ associated to the functor $\overline{p}_{M}:{\rm Tw}(M)\to \sh_{\infty}(\mf/M, \mathscr{S}{\rm p})$, must satisfy descent on $M$. Thus, we have the following.

\begin{proposition}
The presheaf $M\mapsto \big(p_{M}:\lambda_{\widehat{\R}}(M)\longrightarrow {\rm Tw}_{\widehat{\R}}(M)\big)$ satisfies descent.
\end{proposition}

We now would like to use the canonical bundle over the stack of twists to understand how twisted differential cohomology theories are pieced together from local data. With the powerful machinery developed by Lurie, this is now fairly systematic to construct.

\begin{theorem}
[Local-to-global construction of twisted differential cohomology]
\label{Thm-equiv}
Let $M$ be a smooth manifold and let $\hat{\tau}:M\to {\rm Tw}_{\widehat{\R}}$ be a twist for a differential refinement $\widehat{\R}$ of a ring spectrum $\R$ such that the underlying topological twist $\tau$ lands in the identity component of ${\rm Pic}^{\rm top}_{\R}$. Then the following sets are in bijective correspondence:
\begin{enumerate}[label=(\roman*)]

\vspace{-2mm}
\item Equivalence classes of differential twists $\widehat{\R}_{\hat{\tau}}$ for a refinement $\widehat{\R}$.

\vspace{-2mm}
\item Equivalence classes of pullback bundles $E\to M$ by classifying maps $\hat{\tau}:M\to {\rm Tw}_{\widehat{\R}}$
$$
\xymatrix@R=1.5em{
E\ar[rr]\ar[d] && \lambda_{\widehat{\R}}\ar[d]
\\
M\ar[rr]^{\hat{\tau}} && {\rm Tw}_{\widehat{\R}}\;.
}
$$

\vspace{-2mm}
\item Equivalence classes of colimits of the form
$$
\left\{\vcenter{\xymatrix{ E\ar[d]
\\
M
}}\right\}\simeq 
{\rm colim}\left\{\vcenter{\xymatrix@R=1.4em{
\hdots \ar@<.15cm>[r] \ar[r] \ar@<-.15cm>[r] &\coprod_{\alpha\beta}\widehat{\R}\times U_{\alpha\beta}\ar[d]\ar@<.1cm>[r] \ar@<-.1cm>[r]& \coprod_{\alpha}\widehat{\R}\times U_{\alpha}\ar[d]
\\
\hdots \ar@<.15cm>[r] \ar[r] \ar@<-.15cm>[r] &\coprod_{\alpha\beta}U_{\alpha\beta} \ar@<.1cm>[r] \ar@<-.1cm>[r]& \coprod_{\alpha}U_{\alpha}
}}\right\}\;,
$$
where the simplicial maps in the top row are completely determined by a commutative diagram

\vspace{-4mm}
\(\label{exptw1}
\xymatrix@R=1.5em{
\vdots & \vdots
\\
\Delta[1]\ar@<.15cm>[u] \ar[u] \ar@<-.15cm>[u]  \ar[r] & \prod_{\alpha\beta}{\rm Tw}_{\widehat{\R}}(U_{\alpha\beta}) \ar@<.15cm>[u] \ar[u] \ar@<-.15cm>[u] 
\\
\Delta[0] \ar@<.1cm>[u] \ar@<-.1cm>[u]\ar[r] & \prod_{\alpha}{\rm Tw}_{\widehat{\R}}(U_{\alpha})\ar@<.1cm>[u] \ar@<-.1cm>[u]
\;,
}
\)
with $\Delta[0]\longrightarrow \prod_{\alpha}{\rm Tw}_{\widehat{\R}}(U_{\alpha})$ picking out the trivial twist $\widehat{\R}$ on each $U_{\alpha}$.
\end{enumerate}
\end{theorem}
\theproof
The bijection between (i) and (ii) is follows simply Yoneda Lemma
$
\map(M,{\rm Tw}_{\widehat{\R}})\simeq {\rm Tw}_{\widehat{\R}}(M)
$,
along with the fact that the objects of ${\rm Tw}_{\widehat{\R}}(M)$ are by definition precisely the twists $\widehat{\R}_{\hat{\tau}}$ of $\widehat{\R}$. By definition, the global sections of $E$ fit into a diagram \footnote{In the diagram, $\map_{M}(M,E)$ denotes the mapping spectrum in the $\infty$-slice over $M$.}
$$
\xymatrix@R=1.5em{
\Gamma(M;E):=\map_{M}(M,E)\ar[rr]\ar[d] && \lambda_{\widehat{\R}}(M)\ar[d]
\\
\ast \ar[rr]^-{\hat{\tau}} && {\rm Tw}_{\widehat{\R}}(M)
}
$$
and, by the $\infty$-Grothendiek construction, the fiber of $\lambda_{\widehat{\R}}(M)\to  {\rm Tw}_{\widehat{\R}}(M)$ is $\widehat{\R}_{\hat{\tau}}(M)\simeq \Gamma(M;E)$.

For the bijection between (ii) and (iii) recall that, by Proposition \ref{prop-triviali}, for each $x\in M$, 
there is $x\in U\subset M$ such that each differential twist $\widehat{\R}_{\hat{\tau}}$ 
on $U$ is equivalent to the trivial twist $\widehat{\R}$. Thus for each such $U$ we can choose a representative of the class of $\hat{\tau}$ so that the corresponding element in ${\rm Tw}_{\widehat{\R}}(U)$ is $\widehat{\R}$.  Now let $\{U_{\alpha}\}$ be a good open covering 
of such open sets. Consider the augmented colimiting diagram in $\sh_{\infty}(\mf)$
$$
\xymatrix{
\hdots \ar@<.15cm>[r] \ar[r] \ar@<-.15cm>[r] &\coprod_{\alpha\beta}U_{\alpha\beta} \ar@<.1cm>[r] \ar@<-.1cm>[r]& \coprod_{\alpha}U_{\alpha} \ar[r] &M\ar[r]^-{\hat{\tau}} & {\rm Tw}_{\widehat{\R}}
}\;.
$$
By the preceding discussion this gives rise to a colimiting diagram in the tangent topos $T(\sh_{\infty}(\mf))$ via the global sections functor. Taking iterated pullbacks by the universal bundle $\lambda_{\widehat{\R}}\to {\rm Tw}_{\widehat{\R}}$, we get a commutative diagram 
\(\label{dsc ubun}
\xymatrix@R=1.5em{
\hdots \ar@<.15cm>[r] \ar[r] \ar@<-.15cm>[r] &\coprod_{\alpha\beta}\widehat{\R}\times U_{\alpha\beta}\ar[d]\ar@<.1cm>[r] \ar@<-.1cm>[r]& \coprod_{\alpha}\widehat{\R}\times U_{\alpha}\ar[d] \ar[r] & E\ar[d] \ar[r] & \lambda_{\widehat{\R}}\ar[d]
\\
\hdots \ar@<.15cm>[r] \ar[r] \ar@<-.15cm>[r] &\coprod_{\alpha\beta}U_{\alpha\beta} \ar@<.1cm>[r] \ar@<-.1cm>[r]& \coprod_{\alpha}U_{\alpha} \ar[r] &M\ar[r]^-{\hat{\tau}} & {\rm Tw}_{\widehat{\R}}\;.
}
\)
The tangent topos itself is always an $\infty$-topos (see \cite[section 6.1.1]{Lur}). By the Pasting Lemma for pullbacks, all the square diagrams are Cartesian. By the Descent Axiom for a topos, the top diagram must be colimiting, since the bottom diagram is colimiting. Finally, using descent for the stack ${\rm Tw}_{\widehat{\R}}$, we have an equivalence of $\infty$ groupoids
$$
\lim\left\{\vcenter{
\xymatrix{
\hdots & \ar@<.15cm>[l] \ar[l] \ar@<-.15cm>[l] \prod_{\alpha}{\rm Tw}_{\widehat{\R}}(U_{\alpha\beta}) & \ar@<.1cm>[l] \ar@<-.1cm>[l] \coprod_{\alpha}{\rm Tw}_{\widehat{\R}}(U_{\alpha})
}}
\right\}\simeq {\rm Tw}_{\widehat{\R}}(M)\;.
$$
The limit on the left can be calculated by the local formula for the homotopy limit. This gives a bijective correspondence between the classes of $\hat{\tau}$ and the commutative diagrams of the form \eqref{exptw1}, taken up to equivalence.
\endofproof

\subsection{Twisted differential theories as flat bundles of spectra}
\label{Sec-flat}

Theorem \ref{Thm-equiv} gives us several different, but equivalent, ways to think about twisted differential cohomology theories. When taken altogether, one can find a long list of analogies between these objects and the notion of a line bundle equipped with flat connection. Ultimately, we would like to say that a twisted differential theory is a bundle of spectra equipped with flat connection. However, the notion of connection on a bundle of spectra has not yet been defined, nor is it completely clear how natural such a definition would be. For this reason, we will first motivate the notion via a list of analogies. We will then define the connections rigorously and show that indeed twisted differential theories come equipped with a canonical connection as part of the data.

\medskip
\noindent {\bf I. \underline{Local triviality of line bundles}}

\medskip
Consider a locally trivial, real line bundle over a smooth manifold $p:\eta\to M$. The local triviality means that if we fix a good open cover $\{U_{\alpha}\}$ on $M$ and we are given the choice of transition functions $g_{\alpha\beta}$ on intersections $U_{\alpha}\cap U_{\beta}$, then we can piece together the total space of the bundle via the local trivializations. Moreover, a choice of {\v C}ech cocycle $g_{\alpha\beta}$ on intersections can be described in the language of smooth stacks as a commutative diagram 

\vspace{-5mm}
$$
\xymatrix@R=1.5em{
\vdots & \vdots
\\
\Delta[1]\ar@<.15cm>[u] \ar[u] \ar@<-.15cm>[u]  \ar[r]^-{g_{\alpha\beta}} 
& \prod_{\alpha\beta}\BB \RR^{\times}(U_{\alpha\beta}) \ar@<.15cm>[u] 
\ar[u] \ar@<-.15cm>[u] 
\\
\Delta[0] \ar@<.1cm>[u] \ar@<-.1cm>[u]\ar[r] & \prod_{\alpha}\BB \RR^{\times}(U_{\alpha})\ar@<.1cm>[u] \ar@<-.1cm>[u]
\;,
}
$$
with the vertical maps on the right given by restriction and the maps on the left are the face inclusions. By descent, this data is equivalently the data of a map
(see \cite{FSSt} for similar constructions) 
$
g:M\longrightarrow \BB\RR^{\times}
$.
We can construct a universal line bundle over $\BB\RR^{\times}$ via an action of $\RR^{\times}$ on $\RR$. The quotient by this action leads to a stack $\RR/\!/\RR^{\times}$ and we have a canonical map $\RR/\!/\RR^{\times}\to \BB\RR^{\times}\;,$ which projects out $\RR$. The statement of local triviality is then translated into a statement about descent by considering the diagram
\(\label{lbun desc}
\xymatrix@R=1.5em{
\coprod_{\alpha\beta}\RR \times U_{\alpha\beta}\ar[d]\ar@<.1cm>[r] \ar@<-.1cm>[r]& \coprod_{\alpha}\RR \times U_{\alpha}\ar[d] \ar[r] & E\ar[d] \ar[r] & \RR/\!/\RR^{\times}\ar[d]
\\
\coprod_{\alpha\beta}U_{\alpha\beta} \ar@<.1cm>[r] \ar@<-.1cm>[r]& \coprod_{\alpha}U_{\alpha} \ar[r] &M\ar[r]^-{g} & \BB\RR^{\times}
\;.
}
\)
If $E\to M$ is the pullback of the universal bundle map by $g$, and if all the other squares are Cartesian, then descent says that the top diagram is colimiting. Conversely, if all squares but possibly the first are Cartesian and the top diagram is colimiting, then the first square is Cartesian. This means that the map $g$ gives precisely the data necessary to construct the total space of the bundle $E$ over $M$. 

\medskip
\medskip
\noindent {\bf  I${}^\prime$. \underline{Local triviality of bundles of spectra}}
\medskip

Comparing diagram \eqref{lbun desc} with diagram \eqref{dsc ubun}, we see a close analogy between the concept of a locally trivial line bundle and that of
 bundle of spectra $E\to M$ classified by a twist $\hat{\tau}$. In the case of a twisted differential theory, the cocycles representing the transition functions are replaced by maps to the stack of twists. These maps encode how to piece together the total space of the bundle via automorphisms of the underlying theory. Given that we have local triviality available, almost all the constructions one can do for vector bundles go through the same way for bundles of spectra, where the basic operations are now operations on spectra, rather than on vector spaces. 
 
\begin{definition}
[Smash product of bundles of spectra]
Let $E\to M$ and $F\to M$ be two bundles of spectra over a smooth manifold $M$ with fibers $\R$ and $\mathscr{K}$, respectively. We define the \emph{smash product} of $E$ and $F$ as the bundle $E\wedge F\to M$ with fiber $\R\wedge \mathscr{K}$, given by the colimit 
$$
\xymatrix@R=1.5em{
\hdots \ar@<.15cm>[r] \ar[r] \ar@<-.15cm>[r] &\coprod_{\alpha\beta}\R\wedge \mathscr{K} \times U_{\alpha\beta}\ar[d]\ar@<.1cm>[r] \ar@<-.1cm>[r]& \coprod_{\alpha}\R\wedge \mathscr{K}\times U_{\alpha}\ar[d] \ar[r] & E\wedge F\ar[d]
\\
\hdots \ar@<.15cm>[r] \ar[r] \ar@<-.15cm>[r] &\coprod_{\alpha\beta}U_{\alpha\beta} \ar@<.1cm>[r] \ar@<-.1cm>[r]& \coprod_{\alpha}U_{\alpha} \ar[r] & M\;,
}
$$
where the maps are induced by the corresponding maps for $\R$ and $\mathscr{K}$ via the smash product operation $$
\wedge:\sh_{\infty}(\mf;\mathscr{S}{\rm p})\times \sh_{\infty}(\mf;\mathscr{S}{\rm p})
\longrightarrow \sh_{\infty}(\mf;\mathscr{S}{\rm p})\;.
$$
\end{definition}

\medskip

\noindent {\bf II. \underline{Flat connections on a line bundle}}

\medskip
Now we would like to add flat connections to the picture. Let $\nabla$ be a flat connection on the line bundle $L\to M$. A connection is an operator on the sheaf of local sections
$$
\nabla:\Gamma(-;L)\longrightarrow \Gamma(-;L\otimes T^*M)\;,
$$
which satisfies the Leibniz rule with respect to the module action of $C^{\infty}(-;\RR)$ on $\Gamma(-; L)$. If the connection is flat, then $\nabla$ gives rise to a complex
$$
{\cal L}_L:=\big(
\xymatrix{
0 \ar[r] & \Gamma(-;L)\ar[r]^-{\nabla} & \Gamma(-; L\otimes T^*M)\ar[r]^-{\nabla} & \Gamma(-; L\otimes \Exterior^2(TM))\ar[r] & \hdots
}
\big)\;.
$$
With the discussion right after Definition \ref{def-dRA} in mind, we have the following.

\begin{lemma} 
The complex ${\cal L}_L$ is an invertible module over the de Rham complex. 
\end{lemma}
\theproof
To see this, we consider the de Rham complex as an unbounded sheaf of chain complexes concentrated in nonpositive degrees. Then the tensor product 
${\cal L}_L\otimes \Omega^*$ is computed as the complex with degree $-n$ given by 
$$
\bigoplus_{p+q=-n}\Gamma
\big(-; L\otimes \Exterior^p(TM)\big)\otimes
 \Gamma \big(-;\Exterior^q(TM)\big)\;.
$$
The wedge product of forms then gives a map 
$$
\mu:\bigoplus_{p+q=n}\Gamma\big(-; L\otimes \Exterior^p(TM)\big)
\otimes \Gamma\big(-;\Exterior^q(TM)\big)\longrightarrow
\Gamma\big(-; L\otimes \Exterior^n(TM)\big)\;.
$$
To get a map at the level of complexes, we need to check the commutativity with the differentials. But this follows immediately from the Leibniz rule:
\bea
\mu \big((\nabla\otimes {\rm id}\pm {\rm id}\otimes d)(s\otimes \omega\otimes \nu)\big) 
&=& \nabla(s\otimes \omega)\wedge \nu\pm (s\otimes \omega)\wedge d\nu
\\
&=& \nabla(s\otimes \omega\wedge \nu)
\\
&=& \nabla \big(\mu(s\otimes \omega\otimes \nu)\big)\;.
\eea
Thus we have a module action of $\Omega^*$ on the complex ${\cal L}_L$. This complex is invertible with inverse given by the complex ${\cal L}^{-1}_L$, defined to be degreewise identical to ${\cal L}_L$ but equipped with the differential $-\nabla$. 
\endofproof

\medskip
\noindent {\bf II${}^\prime$. \underline{Flat connections on a bundle of spectra}}
\medskip

Continuing with the analogy for twisted differential theories, we see that we can think of the K-flat, invertible module ${\cal L}$ (see the description right after Def.
\ref{def-dRA}) as the sheaf of local sections of the bundle of spectra 
$p: E\wedge H(\Omega^*) \to M$, where $E \to M$ is the underlying bundle of spectra corresponding to the topological twist.  
Being a complex, ${\cal L}$ comes equipped with a differential, which we denote by 
$\nabla$. We will continue to make a distinction between a bundle of spectra 
$E_{\tau}$ and the corresponding sheaf of spectra $\R_\tau$ given by its local sections. 

\medskip
Before defining connections on bundles of spectra, we will need some set up. Let $\Exterior^{\bullet}(T^*M)$ denote the exterior power bundle, viewed as a graded vector bundle over $M$. The sheaf of sections of this graded bundle gives the graded algebra of differential forms on $M$. Equipping this graded algebra with the exterior derivative gives a DGA, and an application of the Eilenberg-Maclane functor gives a corresponding sheaf of spectra on $M$. To simplify notation and highlight our analogy with bundles, we will denote the corresponding bundle of spectra by the same symbol $\Exterior^{\bullet}(T^*M)$. 

\medskip
Now given a bundle of spectra $E\to M$, the sheaf of sections of the smash product $E\wedge \Exterior^{\bullet}(T^*M)$ admits the structure of a sheaf of $H\RR$-module spectrum, inherited from $\Exterior^{\bullet}(T^*M)$. By Shipley's theorem \cite{Sh}, there is corresponding sheaf of DGA's ${\cal L}$, which is unique up to equivalence such that 
$$H({\cal L})\simeq \Gamma(-;E\wedge \Exterior^{\bullet}(T^*M))\;,$$
so that the sections of $\Gamma(-;E\wedge \Exterior^{\bullet}(T^*M))$ admit the structure of a sheaf of DGA's. This structure is not unique, but is unique up to a contractible choice. From the definition of $\Exterior^{\bullet}(T^*M)$, we also see that ${\cal L}$ can be chosen so that it admits the structure of a module over the de Rham DGA $\Omega^*$. This now allows us to use the familiar structure of DGA's to construct connections. 

\begin{definition}
[Connection on a bundle of spectra]
Fix a sheaf of DGA's modeling the sections of the bundle $E\wedge \Exterior^{\bullet}(T^*M)$. A \emph{connection} $\nabla$ on a bundle of spectra $E \to M$ is a morphism of sheaves of DGA's
$$
\nabla: \Gamma(-;E\wedge \Exterior^{\bullet}(T^*M))\longrightarrow
 \Gamma(-;E\wedge\Exterior^{\bullet}(T^*M)[1])\;,
$$
satisfying the Leibniz rule
$$
\nabla(\omega \cdot s)=d(\omega)\cdot s + (-1)^{{\rm deg}\;\omega}
\; \omega\cdot \nabla(s)\;,
$$
where ``$\cdot$" denotes the module action by differential forms. Here $\Exterior^{\bullet}(T^*M)[1]$ denotes the shift of the corresponding bundle of spectra.
\end{definition}

The following proposition shows that every differential twist gives a flat connection on the bundle of spectra corresponding to its underlying topological twist.
\begin{proposition}
Every differential twist $\widehat{\R}_{\hat{\tau}}=(\R_{\tau},t,{\cal L})$ defines a flat connection on the bundle of spectra $E_{\tau} \to M$, classified by the topological twist $\tau$. 
\end{proposition}
\theproof
By definition of the differential twist, we have an equivalence
$$t:H(\L)\simeq \R_{\tau}\wedge H\underline{\RR}\;.$$
By the Poincar\'e Lemma, we have an equivalence of sheaves of spectra $H\underline{\RR}\simeq H(\Omega^*)$ so that $\R_{\tau}\wedge H\underline{\RR}\simeq \R_{\tau}\wedge H(\Omega^*)$. Composing equivalences, we see that $\L$ indeed gives the structure of a sheaf of DGA's to the sections of the bundle $E\wedge \Exterior^{\bullet}(T^*M)\to M$. Denote the differential on $\L$ by $\nabla$. Then the shift $\L[1]$ represents the sections of the shift $E\wedge \Exterior^{\bullet}(T^*M)[1]$ and, since $\nabla^2=0$, we see that indeed $\nabla$ defines a flat connection on the bundle.
\endofproof

\medskip
\medskip
\noindent {\bf III. \underline{Reducing the structure group of a flat line bundle}}

\medskip
Finally, to complete the analogy, we need to discuss how the Riemann-Hilbert correspondence enters the picture. For the line bundle $p:L\to M$, if we take local parallel sections, we get a preferred choice of local trivialization of the bundle. To show this, fix a good open cover $\{U_{\alpha}\}$ of $M$ and let $s_{\alpha}$ be local nonvanishing sections of $E$ satisfying $\nabla(s_{\alpha})=0$. Let $\phi_{\alpha}(x,v)=(x,s_{\alpha}(x))$ be the corresponding trivialization on each $U_{\alpha}$. The transition functions $g_{\alpha\beta}$ written in these local trivializations are then constant, since $\nabla(s_{\alpha})=0$. At the level of smooth stacks, this means that the map $g:M\to \BB\RR^{\times}$ can be chosen to factorize through the stack $\BB\underline{\RR}^{\times}$, where $\underline{\RR}^{\times}$ denotes the constant sheaf on $\RR^{\times}$, viewed as a set via the discrete topology. In these local coordinates, the flat connection trivializes. This story can be summarized succinctly in the category of smooth stacks as follows.
\begin{lemma}
We have a homotopy commutative diagram
$$
\xymatrix@R=1.5em{
 M\ar[d]\ar[rr] &&  \flat {\rm Line}_{\nabla}\ar[dll]^-{\rm RH}
\\
\BB\underline{\RR}^{\times} 
}\;,
$$
where the stack $\flat {\rm Line}_{\nabla}$ is the stack of flat line bundles, given locally by Dold-Kan image of the sheaf of complexes
$$
\big(
\xymatrix{
\hdots \ar[r] & 0\ar[r] & C^{\infty}(-;\RR^{\times})\ar[r]^-{d} & \Omega^1_{\rm cl} \ar[r] & 0\ar[r] & \hdots
}\big)\;.
$$
The diagonal map  ${\rm RH}$  is an equivalence,
given by the Riemann-Hilbert correspondence, which associates every flat bundle to its corresponding local system. 
\end{lemma}

More generally, we can use the affine structure for flat connections to add globally defined closed 1-forms to the flat connection and the resulting connection is trivially compatible with the discrete structure. We can also ask for further compatibility on $M$ as follows. Consider the inclusion 
$\ZZ/2\into \RR^{\times}$ as the units of $\ZZ$. If we ask for the flat connections to have monodromy representation factoring through $\ZZ/2$, then we are led to a homotopy commutative diagram 
$$
\xymatrix@R=1.5em{
 M\ar[d]\ar[rr] &&  \flat {\rm Line}_{\nabla}\ar[d]^-{\rm RH}
\\
\BB\ZZ/2 \ar[rr] && \BB\underline{\RR}^{\times}
\;.
}
$$
If we further ask that the chosen connection is globally defined and compatible with the $\ZZ/2$-structure, then we are asking for a homotopy commutative diagram 
$$
\xymatrix@R=1.5em{
 M\ar[d]\ar[rr] &&  \Omega^1_{\rm cl} \ar[d]
\\
\BB\ZZ/2 \ar[rr] && \BB\underline{\RR}^{\times}\simeq \flat{\rm Line}_{\nabla}
\;.
}
$$
Then, by the universal property of the homotopy pullback, we have an induced map 
\footnote{The superscript $h$ is for ``homotopy", and this should not be
confused with the class of twists $h$.}
$$
M\longrightarrow \BB(\ZZ/2)_{\nabla}:=
\BB\ZZ/2\times^h_{{}_{\BB\underline{\RR}^{\times}}}\Omega^1_{\rm cl}\;.
$$
\begin{remark} 
The stack $\BB(\ZZ/2)_{\nabla}$ is thus the moduli stack of line bundles equipped with globally defined flat connection, with monodromy factoring through $\ZZ/2$. Equivalently, it is the moduli stack of those (globally) flat bundles whose transition functions (when written in the trivializations provided by parallel sections) take values in $\ZZ/2$. 
\end{remark} 

\medskip
\noindent {\bf III${}^{\prime}$. \underline{Reducing the structure group of 
a flat bundle of spectra}}

\medskip
The analogue of the Riemann-Hilbert correspondence for twisted differential theories is essentially the twisted de Rham theorem, which gives an equivalence of stacks
$$
{\rm Pic}^{\rm dR}_{\R}\simeq{\rm Pic}^{\rm top}_{\R\wedge H\RR}\;.
$$
Geometrically, we begin by solving the equation $\nabla(s)=0$ locally on $M$, where $\nabla$ is the differential on the K-flat invertible module $\L$. In general, we need to know that such solutions exist locally and, for this, we need to better understand the structure of $\L$.

\begin{proposition}
[Structure of invertible modules over ${\cal R}$]
\label{Prop-invert}
If $\L$ is an invertible module over ${\cal R}$, then $\L$ is locally isomorphic to a product of (possibly shifted) finitely many copies of ${\cal R}$.
\end{proposition}
\theproof
Let ${\cal K}$ invert $\L$, i.e. be such that $\L\otimes_{{}_{{\cal R}}}{\cal K}\cong {\cal R}$ is an isomorphism of sheaves of DG ${\cal R}$-modules. Then there is a point $x\in U\subset M$ such that, in that neighborhood, 
$\L(U)\otimes_{{}_{{\cal R}(U)}}{\cal K}(U)\cong {\cal R}(U)$.
Degree-wise this reduces to a collection of isomorphisms
$$
\xymatrix{
\phi_{n}:\bigoplus_{p+q=n}\L(U)_p\otimes_{{}_{{\cal R}(U)_0}}
{\cal K}(U)_{q}/\sim 
\; \; \ar[r]^-\cong &
 {\cal R}_n(U)\;,
 }
$$
where $\sim$ denotes the submodule generated by the relation $(rl,k)=(l,rk)$, with $\vert l\vert+\vert r \vert+\vert k \vert=n$. Let $1\in {\cal R}(U)_0$ denote the identity. Then, for $n=0$, there is a finite sum 
$\sum_{i,p}l_{p,i}\otimes k_{-p,i}$ (with the $k_{p, i}$'s and $l_{p, i}$'s 
independent in $\L_{p}$ and ${\cal K}_p$, respectively) that is mapped as $\phi_{0}(\sum_{i,p}l_{p,i}\otimes k_{-p,i})=1$. Since $d(1)=0$ and $\phi$ is 
a chain map, we must have
$$
0=d(1)=d\phi_{0}\big(\sum_{i,p}l_{p,i}\otimes k_{-p,i}\big)=\phi_{-1}\Big(d\big(\sum_{i,p}l_{p,i}\otimes k_{-p,i}\big)\Big)\;.
$$
Since $\phi$ is an isomorphism, this implies that 
$
0=d\big(\sum_{i,p}l_{p,i}\otimes k_{-p,i}\big)=\sum_{p,i}
dl_{p,i}\otimes k_{p,i}\pm l_{p,i}\otimes dk_{-p,i}
$.
Since the $k_{-p,i}$'s and $l_{p,i}$'s were chosen to be independent, this forces $dl_{p,i}=0$ and $dk_{-p,i}=0$ for each $p$ and $i$. Define the complex
${\cal L}^{\prime}(U)$ as the subcomplex of ${\cal L}(U)$ which is freely generated by $l_{p,i}$ as a graded module over ${\cal R}(U)$. Its underlying chain complex can be identified with the complex
$$
 \bigoplus_{p,i}{\cal R}(U)[p]\cong \bigoplus_{p}~\langle l_{p,i}\rangle\;,
$$
where $\langle l_{p,i}\rangle$ is the freely generated ${\cal R}(U)$-module on generators $l_{p,i}$ (as $i$ varies) shifted up by $p$. The differential $d$ on this complex is inherited from ${\cal L}(U)$ and, therefore, satisfies $d^2=0$. Now the tensor product of the inclusion $i:{\cal L}^{\prime}(U)\into {\cal L}(U)$ (which is a degreewise monomorphism) with ${\cal K}(U)$ is also a degree-wise epimorphism. Indeed, the element $\sum_{i,p}l_{p,i}\otimes k_{-p,i}$ generates the tensor product as an ${\cal R}(U)$-module. By invertibility, we must have a commutative diagram
$$
\xymatrix{
\L^{\prime}(U)\ar[d]^-{\cong}\ar[rr]^{i} && \L(U)\ar[d]^-{\cong}
\\
\L^{\prime}(U)\otimes_{{}_{{\cal R}(U)}}{\cal K}(U)
\otimes_{{}_{{\cal R}(U)}} \L(U)\ar[rr] &&  {\cal L}(U)
\otimes_{{}_{{\cal R}(U)}} {\cal K}(U)\otimes_{{}_{{\cal R}(U)}}
 \L(U)
\;,
}
$$
with the bottom map a degree-wise epimorphism via right exactness of the tensor product. Thus $i$ is an isomorphism. 
\endofproof

Now focusing on the special case where ${\cal R}(U)=\Omega^*(U;A)$, we have the following.
\begin{proposition}
\label{Prop-invR}
Every invertible module ${\cal L}$ over $\Omega^*(-;A)$ is locally isomorphic to a shift $\Omega^*(-;A)[p]$.
\end{proposition}
\theproof
We need to show that, in this case, we can choose the $l_{p,i}$ in the proof of Proposition \ref{Prop-invert} to have a single generator in each degree (i.e. the indexing set for $i$ contains one element). To this end, recall that (locally) the ring $C^{\infty}(-;\RR)$ is freely generated by a single nonvanishing section. By Proposition \ref{Prop-invert}, $\L$ and ${\cal K}$ locally look like a finite product of shifts of $\Omega^*(-;A)$. After a careful consideration of the grading for ${\cal L}$, ${\cal K}$ and $\Omega^*(-;A)$, we find that there is $x\in U\subset M$ and free $C^{\infty}(U;\RR)$-submodules 
${\cal L}^{\prime}(U)_p\cong \bigoplus_{i} C^{\infty}(U;\RR)$ and ${\cal K}^{\prime}(U)_{-p}\cong \bigoplus_{j} C^{\infty}(U;\RR)$ of ${\cal L}(U)_{p}$ 
and ${\cal K}(U)_{-p}$, respectively, such that 
$$
\bigoplus_{p}{\cal L}^{\prime}(U)_p\otimes_{{}_{C^{\infty}(U;\RR)}}
 {\cal K}^{\prime}(U)_{-p}\cong  C^{\infty}(U;\RR)\;.
$$
i.e. the only elements which can get mapped to smooth functions in $\Omega^*(-;A)$ must come from the given submodules. Using basic properties of the tensor product, along with the fact that the isomorphism is that of $C^{\infty}(U;\RR)$-modules, implies that the indexing sets for $i$, $j$, and $p$, must contain a single element.
\endofproof

\begin{remark}
Notice that Proposition \ref{Prop-invR} also implies that the condition of 
``weakly locally constant"  can be dropped from the definition of the stack ${\rm Pic}^{\rm form}_{\R}$ in Section \ref{Sec-rev}, as $\Omega^*(-;A)$ and its shifts 
are weakly locally constant by de Rham's Theorem.
\end{remark}

As a consequence of Proposition \ref{Prop-invR}, we can immediately 
solve our differential equation.

\begin{corollary}
The equation $\nabla(s^{p}_{\alpha})=0$ admits a solution which is provided by the isomorphism in Proposition \ref{Prop-invR}. Moreover, $s_{\alpha}$ gives rise to a local isomorphism via the module action
$$
s_{\alpha}\wedge (-):\Omega_{U_{\alpha}}^*(-;A)[p]\overset{\cong}{\longrightarrow} \L_{U_{\alpha}}\;.
$$
\end{corollary} 
\theproof
In the proof of Proposition \ref{Prop-invR}, take $s_{\alpha}^{p}=l_{p}$.
\endofproof

It is natural to ask about the resulting structure globally. For that, if we compare 
the isomorphisms $s_{\alpha}$ on intersections, we get an automorphism
$$
\xymatrix@R=1.5em{
\Omega_{U_{\alpha\beta}}^*(-;A)\ar[rr]^-{g_{\alpha\beta}}\ar[dr]_{s_{\alpha}\wedge} && \Omega_{U_{\alpha\beta}}^*(-;A)\;.
\\
& \L_{U_{\alpha\beta}}\ar[ur]_{s_{\beta}^{-1}\wedge}&
}
$$

\vspace{-2mm}
\noindent The inverse map $s_{\beta}^{-1}\wedge$ takes a section of the form 
$s_{\beta}\wedge \omega$ 
to its coefficient $\omega$. Thus, the automorphism $g_{\alpha\beta}$ is represented as a wedge product with a form. Moreover, since both $s_\alpha$ and 
$s_\beta^{-1}$ are chain maps, this form is closed
$$
dg_{\alpha\beta}=d(s_{\beta}^{-1}s_{\alpha})=s_{\beta}^{-1}(\nabla(s_{\alpha}))=0\;.
$$
By the Poincar\'e Lemma, it is also exact. This results in a degree one automorphism 
$s_{\alpha\beta}$, represented by a wedge product with a form
$$
d(s_{\alpha\beta}\wedge \omega)\pm s_{\alpha\beta}\wedge d\omega
=g_{\alpha\beta}\wedge \omega\;.
$$
Continuing to higher intersections, we get a whole hierarchy of automorphisms on various intersections which are compatible in a prescribed way with the automorphisms one step down. The compatibility is precisely that of an element in the total complex of the {\v C}ech-de Rham 
double complex with values in ${\cal L}$. More explicitly, on the $n$-fold intersecting patch $U_{\alpha_0\hdots \alpha_{n}}$, the sections $s_{\alpha_0\hdots \alpha_n}$ is closed under the total differential of the double complex, i.e.,   
$$
\nabla(s_{\alpha_0\hdots \alpha_n})=\pm \delta s_{\alpha_0\hdots \alpha_{n-1}}\;.
$$
We now use these sections to consider lifting through diagram \eqref{bun twpull}.

\begin{theorem}
[Twisted differential spectra from local data]
\label{thm-sections}
Let $\widehat{\R}=(\R,{\rm ch},A)$ be a differential cohomology theory and let ${\cal L}$ be a K-flat, invertible module over $\Omega^*(-;A)$ with $s=(s_{\alpha},s_{\alpha\beta},\hdots )$ the corresponding cocycle with values in the {\v C}ech total complex of ${\cal L}$. Assume that the isomorphism $s_{\alpha}$ identifies ${\cal L}$ with the unshifted complex $\Omega^*(-;A)$ locally. Then $s$ determines a commutative diagram

\vspace{-5mm}
$$
\xymatrix@R=1.5em{
\vdots & \vdots
\\
\Delta[1]\ar@<.15cm>[u] \ar[u] \ar@<-.15cm>[u]  \ar[r]^-{\eta_{\alpha\beta}} & \prod_{\alpha\beta}B{\rm GL}_1(\R\wedge H\RR)\ar@<.15cm>[u] \ar[u] \ar@<-.15cm>[u] 
\\
\Delta[0] \ar@<.1cm>[u] \ar@<-.1cm>[u]\ar[r]& \prod_{\alpha}B{\rm GL}_1(\R\wedge H\RR) \ar@<.1cm>[u] \ar@<-.1cm>[u]
\;.
}
$$
Moreover, if the $\eta$'s further lift through the canonical map $BGL_1(\R)\to BGL_1(\R\wedge H\RR)$, then $s$ determines a commutative diagram

\vspace{-15mm}
\(\label{diff twcomh}
\xymatrix@R=1.5em{
\vdots & \vdots
\\
\Delta[1]\ar@<.15cm>[u] \ar[u] \ar@<-.15cm>[u]  \ar[r]^-{h_{\alpha\beta}} & \prod_{\alpha\beta}{\rm Tw}_{\widehat{\R}}(U_{\alpha\beta}) \ar@<.15cm>[u] \ar[u] \ar@<-.15cm>[u] 
\\
\Delta[0] \ar@<.1cm>[u] \ar@<-.1cm>[u]\ar[r]^-{h_{\alpha}} & \prod_{\alpha}{\rm Tw}_{\widehat{\R}}(U_{\alpha})\ar@<.1cm>[u] \ar@<-.1cm>[u]
\;,
}
\)
and hence a twisted differential cohomology theory $\widehat{\R}(-; {\hat{\tau}})$, by Theorem \ref{Thm-equiv}.
\end{theorem}
\theproof
For each $\alpha$, the map $s_{\alpha}$ is an isomorphism of $\Omega^*(-;A)$-modules
$$
s_{\alpha}:\Omega_{U_{\alpha}}^*(-;A)
\overset{\cong}{\longrightarrow}{\cal L}_{U_{\alpha}}
$$
over the patch $U_{\alpha}$. Precomposing this map with the canonical inclusion $A\into \Omega^*(-;A)$ yields a morphism $a_{\alpha}$ of $\Omega^*(-;A)$ modules and hence determines an edge in $\prod_{\alpha}{\rm Pic}^{\rm dR}_{\widehat{\R}}(U_{\alpha})$.  Consider the diagram 
determined by the transition forms $g_{\a \b}$
$$
\xymatrix{
\R\wedge H\RR \ar[rr]^-{\rm ch} \ar[d]^-{\eta_{\alpha\beta}} && 
H\big(\Omega_{U_{\alpha\beta}}^*(-;A)\big)\ar[rr]^-{s_{\alpha}}\ar[d]^-{g_{\alpha\beta}} && 
H({\cal L}_{U_{\alpha\beta}})\ar@{=}[d]
\\
\R\wedge H\RR \ar[rr]^-{\rm ch} &&
H\big(\Omega_{U_{\alpha\beta}}^*(-;A)\big)
\ar[rr]^-{s_{\beta}} \ar@/^1pc/[ll] && 
H({\cal L}_{U_{\alpha\beta}})
\;,}
$$
with the bottom curved arrow giving a homotopy inverse for the Chern character map 
and the composition determining the map $\eta_{\alpha\beta}$. 
A choice of differential form $s_{\a \b}$
in $\Omega^*(U_{\alpha\beta};A)$ satisfying $ds_{\alpha\beta}=g_{\alpha\beta}$ determines a homotopy that fills this diagram. Indeed, such a form gives a degree one automorphism
\(\label{cyctr com}
s_{\alpha\beta}\wedge(-) \in Z\Big(\hom_{\mathscr{C}{\rm h}}
\big(\Omega_{U_{\alpha}}^*(-;A),\; \Omega_{U_{\alpha}}^*(-;A)\big)\Big)_1\;,
\)
via the module action, where $Z$ is the functor $\mathscr{C}{\rm h}\to \chp$ which truncates an unbounded complex and puts cycles in degree zero. Since $ds_{\alpha\beta}=g_{\alpha\beta}$, one easily checks that the image of the automorphism $s_{\alpha\beta}\wedge(-)$ under the differential on the complex \eqref{cyctr com} is indeed the automorphism of degree zero given by $g_{\alpha\beta}\wedge (-)$. Applying the Dold-Kan functor to the above positively graded chain complex gives an edge in the automorphism space
${\rm Aut}(\Omega_{U_{\alpha}}^*(-;A),\Omega_{U_{\alpha}}^*(-;A))$. Then pre- and 
post-composing with ${\rm ch}$ and ${\rm ch}^{-1}$ give rise to an edge in 
${\rm Aut}(\R\wedge H\RR,\R\wedge H\RR)$ with boundary $({\rm id},\eta_{\alpha\beta})$.

The data described above is manifestly that of a map
$$
s_{\alpha\beta}:\Delta[1]\times \Delta[1] \; \longrightarrow \;
 \prod_{\alpha\beta}{\rm Pic}^{\rm dR}_{\widehat{\R}}(U_{\alpha\beta})\;,
$$ 
which connects the restrictions of the maps $a_{\alpha}$ and $a_{\beta}$ on intersections. Continuing the process, one identifies the higher forms $s_{\alpha\beta\gamma}$, etc., 
with higher degree automorphisms and gets compatibility with the lower degree 
automorphisms analogously. 
Now, recalling diagram \eqref{bun twpull},
the resulting structure is precisely a commutative diagram 

\vspace{-6mm}
\(\label{dr comtw}
\xymatrix@R=1.5em{
\vdots & \vdots
\\
\Delta[1]\times \Delta[1]\ar@<.15cm>[u] \ar[u] \ar@<-.15cm>[u]  \ar[r]^-{a_{\alpha\beta}} &\prod_{\alpha\beta}{\rm Pic}^{\rm dR}_{\widehat{\R}}(U_{\alpha\beta}) \ar@<.15cm>[u] \ar[u] \ar@<-.15cm>[u] 
\\
\Delta[1]  \times \Delta[0]\ar@<.1cm>[u] \ar@<-.1cm>[u]\ar[r]^-{a_{\alpha}} & \prod_{\alpha}{\rm Pic}^{\rm dR}_{\widehat{\R}}(U_{\alpha})\ar@<.1cm>[u] \ar@<-.1cm>[u]
\;
\;,
}
\)
such that:

\noindent $\bullet$ the evaluation at one of the two endpoints of 
the 1-simplex $\Delta[1]$ factors through the composite map 
$B{\rm GL}_1(\R\wedge H\RR)\into {\rm Pic}^{\rm top}_{\R\wedge H\RR}\to {\rm Pic}^{\rm dR}_{\widehat{\R}}$, hence lifting to a commutative diagram

\vspace{-4mm}
\(\label{bgl ratctw}
\xymatrix@R=1.5em{
\vdots & \vdots
\\
\Delta[1]\ar@<.15cm>[u] \ar[u] \ar@<-.15cm>[u]  \ar[r]^-{\eta_{\alpha\beta}} & \prod_{\alpha\beta}B{\rm GL}_1(\R\wedge H\RR) \ar@<.15cm>[u] \ar[u] \ar@<-.15cm>[u] 
\\
\Delta[0] \ar@<.1cm>[u] \ar@<-.1cm>[u]\ar[r]& \prod_{\alpha}B{\rm GL}_1(\R\wedge H\RR)\ar@<.1cm>[u] \ar@<-.1cm>[u]
\;,
}
\)
$\bullet$ and the evaluation at the other endpoint factors through the map ${\rm Pic}_{\widehat{\R}}^{\rm form}\to {\rm Pic}^{\rm dR}_{\widehat{\R}}$, lifting to a commutative diagram

\vspace{-6mm}
\(\label{form ctw}
\xymatrix@R=1.5em{
\vdots & \vdots
\\
\Delta[1]\ar@<.15cm>[u] \ar[u] \ar@<-.15cm>[u]  \ar[r]^-{s_{\alpha\beta}} & \prod_{\alpha\beta}{\rm Pic}^{\rm form}_{\widehat{\R}}(U_{\alpha\beta}) \ar@<.15cm>[u] \ar[u] \ar@<-.15cm>[u] 
\\
\Delta[0] \ar@<.1cm>[u] \ar@<-.1cm>[u]\ar[r]^-{s_{\alpha}} & \prod_{\alpha}{\rm Pic}^{\rm form}_{\widehat{\R}}(U_{\alpha})\ar@<.1cm>[u] \ar@<-.1cm>[u]
\;.
}
\)
Finally, if the $\eta$'s factor further through the rationalization map 
$B{\rm GL}_1(\R)\to B{\rm GL}_1(\R\wedge H\RR)$ then, by definition of the twisting stack ${\rm Tw}_{\widehat{\R}}$ (i.e. being a pullback), the diagrams \eqref{bgl ratctw} and \eqref{form ctw}  indeed furnish the commutative diagram \eqref{diff twcomh}.
\endofproof

\begin{remark} 
\label{Rem-table}
Now completing the analogy, we see that the stack ${\rm Pic}_{\R}^{\rm top}$, ${\rm Pic}_{\R}^{\rm dR}\simeq {\rm Pic}^{\rm top}_{\R\wedge H\RR}$ and ${\rm Pic}_{\R}^{\rm form}$ (see Diagram \eqref{bun twpull}) can be thought of as the analogues of the stacks $\BB\ZZ/2$, $\BB\underline{\RR}^{\times}\simeq  \flat {\rm Line}_{\nabla}$ and $\Omega^1_{\rm cl}$, respectively
(see also the first table in the Introduction). 
The stack of differential twists ${\rm Tw}_{\widehat{\R}}$ is then analogous to $\BB\ZZ/2_{\nabla}$. 
This comparison indicates that we should think of the twisting stack ${\rm Tw}_{\widehat{\R}}$ as the moduli stack of bundles of spectra $E\to M$, equipped with a globally defined flat connection encoded by the differential on the invertible, K-flat, locally constant module $\L$, whose corresponding local coefficient system is given by the underlying topological twisted theory $\R(-; {\tau})$.
\end{remark}

\begin{remark}
\label{Rem-equi}
\noindent {\bf (i).} Note that the theory of bundles of spectra can be viewed as a generalization of the theory of ordinary line bundles, so that 
the big table in the introduction is not just an analogy. 
 Indeed, the sheaf of sections of a line bundle can be viewed as a sheaf of DGA's, concentrated in degree zero. Application of the Eilenberg-Maclane functor gives a sheaf of ring spectra, which in turn gives rise to a trivial bundle over any fixed manifold $M$. The fiber of this bundle is then easily computed as the Eilenberg-Maclane spectrum $H\RR$. 

\vspace{0.1 cm}

\noindent {\bf (ii)}. In \cite{GS4}, we found that the stack $\BB (\ZZ/2)_{\nabla}$ gives the twists for smooth Deligne cohomology, so that the constructions here indeed represent a generalization to other generalized cohomology theories.

\end{remark} 

\subsection{The general construction of the spectral sequence 
$\widehat{\rm AHSS}_{\hat \tau}$}
\label{Sec-AHSS}

The general machinery established in \cite{GS3} works in the twisted case as well with some modifications. More precisely, let $\widehat{\R}=(\R,c,A)$ be a differential ring spectrum. If we are given a smooth manifold $M$, equipped with a map 
$$
\widehat{\tau}:M\longrightarrow {\rm Tw}_{\widehat{\R}}
$$
to the stack of twists, then $\widehat{\R}_{\hat \tau}=
(\R_{\tau},{\rm ch}_{\tau},{\cal L})$ is a twisted differential spectrum which lives over the manifold $M$; i.e. it is an object in $\sh_{\infty}(M;\Sp)$. Fix a good open cover $\{U_{\alpha}\}$ of the smooth manifold $M$. By the second property in Proposition \ref{prop twdiff}, the restrictions of $\widehat{\R}_{\tau}$ to the various contractible intersections can be identified with the untwisted differential t
heory $\widehat{\R}$.

\medskip This locally trivializing phenomenon is really the key observation in constructing the differential twisted AHSS and is similar to the familiar analogue in spaces. Indeed, the topological twisted AHSS has $E_2$-page which looks like cohomology with coefficients in the underlying \emph{untwisted} theory, and again this amounts to the fact that the twisted $\R$-theory of a contractible space reduces to untwisted $\R$-theory.

\medskip
We now sketch the filtration in the twisted context. Essentially, this is the same filtration as considered in \cite{GS3}, but where each level of the filtration is equipped with a map to the stack of twists ${\rm Tw}_{\widehat{\R}}$. 

\begin{remark}[The filtration for differential twisted spectra]
\label{Rem-filt}
Let us consider a smooth manifold $M$ and the {\v C}ech nerve of a cover $\{U_{\alpha}\}$ of $M$. The restriction of a twist $\widehat{\tau}:M\to {\rm Tw}_{\widehat{\R}}$ to various intersections of the cover gives a simplicial diagram in the slice $\sh_{\infty}(\mf)/{\rm Tw}_{\widehat{\R}}$. That is, we have the simplicial diagram
$$
\xymatrix@R=1.5em{
\hdots \ar@<.05cm>[r]\ar@/_1pc/[drr] \ar@<-.05cm>[r]\ar@<-.15cm>[r]\ar@<.15cm>[r] & \coprod_{\alpha\beta\gamma}U_{\alpha\beta\gamma}\ar@<.1cm>[r]\ar[r]\ar[dr] \ar@<-.1cm>[r] & \coprod_{\alpha\beta}U_{\alpha\beta}\ar@<.05cm>[r]\ar@<-.05cm>[r] \ar[d] &  \coprod_{\alpha}U_{\alpha}\ar[dl]
\\
&& {\rm Tw}_{\widehat{\R}} &
}\;,
$$ 
where the maps to ${\rm Tw}_{\widehat{\R}}$ are given by the restriction of 
the twisted spectrum $\widehat{\R}_{\hat \tau}$ to various 
intersections. Let $F_p$ denote the skeletal filtration on this simplicial object. In this case, 
the successive quotients $F_p/F_{p-1}$ take the form
$$
F_p/F_{p-1}\simeq \Big({\hat \tau}\vert_{U_{\alpha_0\hdots\alpha_p}}:\Sigma^p\Big(\bigvee_{\alpha_0\hdots\alpha_p}(U_{\alpha_0\hdots\alpha_p})_+\Big)
\longrightarrow {\rm Tw}_{\widehat{\R}}\Big)\;.
$$
Using the properties of $\widehat{\R}_{\hat \tau}=(\R_{\tau},t,{\cal L})$, we see that the connected components of the 
global sections of $\widehat{\R}_{\hat \tau}$  are given by

\vspace{-5mm}
\bea
\widehat{R}^0\Big(\Sigma^p\big(\bigvee_{\alpha_0\hdots\alpha_p}(U_{\alpha_0\hdots\alpha_p})_+\big); {\hat \tau}\Big)
&\cong &
\bigoplus_{\alpha_0\hdots\alpha_p}
\widehat{R}^{-p}( U_{\alpha_0\hdots\alpha_p}; {\hat \tau}) \\
&\cong &
\bigoplus_{\alpha_0\hdots\alpha_p}\widehat{R}^{-p}(U_{\alpha_0\hdots\alpha_p})\;.
\eea
Here the last isomorphism uses local triviality (Proposition 
\ref{prop-triviali}) since $U_{\alpha_0\hdots\alpha_p}$ is contractible as a space. 
\end{remark}

The same arguments used in the ordinary AHSS hold for
 this filtration and we have the following. 

\begin{theorem}
[AHSS for twisted differential spectra]
\label{Thm-ourAHSS}
Let $M$ be a compact manifold and let $\widehat{\R}_{\hat \tau}:M\to {\rm Tw}_{\widehat{\R}}$ 
be a twist for a differential ring spectrum $\widehat{\R}$. Then there is a twisted AHSS type spectral sequence  $\widehat{\rm AHSS}_{\hat \tau}$ with
\begin{eqnarray*}
E^{p,q}_2=H^p\big(M;R^{-q-1}_{U(1)}(\ast)\big)\  &\Rightarrow &
\widehat{R}^{p+q}(M,\hat{\tau}), \; \text{for} ~q<0,
  \\
  E^{p,q}_2=H^p\big(M;R^{-q}(\ast)\big)\ &\Rightarrow & 
   \widehat{R}^{p+q}(M,\hat{\tau}),
\; \text{for} ~q>0,
\\
E^{0,0}_2=H^0(M;\widehat{R}^0(-;{\hat{\tau}}))&\Rightarrow & 
 \widehat{R}^{0}(M,\hat{\tau})\;,
\end{eqnarray*}
where the group $H^0(M;\widehat{R}^0(-;{\hat{\tau}}))$ is the degree zero {\v C}ech cohomology with coefficients in the presheaf $\widehat{R}^0(-;{\hat{\tau}})$, i.e. it is the kernel of the map
\(\label{tw res spR}
\xymatrix{
\prod_{\alpha}\widehat{R}^0(U_{\alpha};{\hat{\tau}})\ar[rr]^-{r_{\alpha\beta}-r_{\beta\alpha}} && \prod_{\alpha\beta}\widehat{R}^0(U_{\alpha\beta};{\hat{\tau}})
}\;,
\)
where $r_{\alpha\beta}$ is the restriction map induced by the inclusion $U_{\alpha\beta}\into U_{\alpha}$ and $r_{\beta\alpha}$ is induced by $U_{\alpha\beta}\into U_{\beta}$. This spectral sequence 
converges to the graded group $E^{p,q}_{\infty}$ associated to the filtration $F_p$ described in Remark \ref{Rem-filt} above.
\end{theorem}

Convergence of spectral sequences is studied in general in \cite{Bo} \cite{Mc}. 
In our case, convergence is understood in the same sense as in the untwisted
differential case \cite{GS3}, similarly to how
in the classical topological case convergence in twisted K-theory 
\cite{Ro1} \cite{Ro} \cite{AS2} has been understood in the same way as for 
untwisted  K-theory \cite{AH}. In the differential case, convergence is guaranteed by virtue of the fact that $M$ is compact and, therefore, admits a finite good open cover. 
This is analogous to the classical topological case, where the underlying space is assumed to be a CW complex of finite dimension. Note, however, that in both
the topological and the differential case this can be extended to certain 
infinite dimensional spaces and manifolds, respectively,  by taking appropriate 
direct limits. However, we will need consider this extension explicitly here. 

\begin{remark}
Notice that the $E^{0,0}$-entry involves the twisted theory as coefficients. This might 
seem to disagree with the local triviality condition for the twisted theory (i.e. one might expect coefficients in the untwisted theory to appear). However, note that moving from the $E_1$-page to the $E_2$-page, we have to find the kernel of the map \eqref{tw res spR} as explained in \cite{GS3}. The isomorphisms which identify the theory locally with the untwisted theory do not commute with the usual restriction maps and so we cannot identify the result with its analogue in the untwisted case. However, we can compute $E^{0,0}_{2}$ alternatively as the kernel
$$
E^{0,0}_2\cong \ker\Big\{\xymatrix{
\prod_{\alpha}\widehat{R}^0(U_{\alpha})\ar[rr]^-{r_{\alpha\beta}-\phi_{\alpha\beta}} && \prod_{\alpha\beta}\widehat{R}^0(U_{\alpha\beta})
}\Big\}\;,
$$
where the maps $\phi_{\alpha\beta}$ are the transition maps induced from the local 
trivializations 
$$
\xymatrix{
\prod_{\alpha\beta}\widehat{R}^0(U_{\alpha\beta})\ar[rr]^{\phi_{\alpha\beta}}\ar[d]_{s_{\alpha}}^-{\cong} && \prod_{\alpha\beta}\widehat{R}^0(U_{\alpha\beta})
\\
\prod_{\alpha\beta}\widehat{R}^0(U_{\alpha\beta},\hat{\tau})\ar@{=}[rr] && \prod_{\alpha\beta}\widehat{R}^0(U_{\alpha\beta},\hat{\tau})
\;.
\ar[u]_{s_{\beta}^{-1}}^-{\cong}
}
$$
\end{remark}

\subsection{Properties of the differentials}
\label{Sec-prop}

In this section, we would like to describe some of the properties which 
the differentials in the $\widehat{\rm AHSS}_{\hat \tau}$ for twisted differential cohomology theories enjoy. Some of these properties involve the module structure of a twisted differential theory over the untwisted theory. We have not yet described the module structure. To this end, we note that for a twisted differential theory $\widehat{\R}_{\tau}=(\R_{\tau},t,{\cal L})$, both $\R_{\tau}$ and ${\cal L}$ are modules over $\R$ and $\Omega^*(-;A)$, respectively. The map $t$ is an equivalence in ${\rm Pic}^{\rm form}_{\widehat{\R}}$ and so is, in particular, a map of module spectra. In the same way one defines a differential refinement of a ring spectrum (see \cite{Bun}\cite{Up}), one sees that the twisted theory 
$\widehat{\R}_{\tau}$ admits module structure over $\widehat{\R}$. At the level of the underlying theory, we have maps
\(\label{mod maprt}
\mu:\widehat{R}^n(M)\times \widehat{R}^{n}(M; \hat{\tau})
\longrightarrow
 \widehat{R}^{n+m}(M; \hat{\tau})\;,
\)
and so at least it makes sense to talk about the behavior of the differentials with respect to the module structure. 

\medskip
Some of the classical properties of the differentials in the spectral sequence for twisted K-theory are discussed 
in \cite{AS1}\cite{AS2}\cite{KS2}.  Many of the analogous properties also hold  for twisted Morava K-theory
\cite{SW}. We will first review these properties and then discuss the differential refinement. Let $\R_{\tau}$ be a twist for a ring spectrum $\R$. In general, we have 

\begin{enumerate}[label=(\roman*)]
\item ({\it Linearity}) Each differential $d_i$ is a $R^*(*)$-module map.

\vspace{-2mm}
\item  ({\it Normalization}) The twisted differential with a zero twist reduces to the untwisted differential (which may be zero).

\vspace{-2mm}
\item ({\it Naturality}) Let $\R_{\tau}$ be a twist of $\R$ on $X$. 
If $f:Y \to X$ is continuous then $f^*\R_\tau$ is a twist of $\R$ on $Y$
and $f$  induces a map of spectral sequences from the twisted $\op{AHSS}$ 
for $\R_{\tau}$ to the twisted AHSS
 for $f^* \R_{\tau}$.  On $E_2$-terms it is induced by $f^*$ in cohomology, and on 
$E_\infty$-terms it is the associated graded map induced by $f^*$ in the
 twisted $\R$ theory $\R_{\tau}$.

\vspace{-2mm}
\item  ({\it Module}) The $\op{AHSS}_\tau$ for $\R_{\tau}$ is a spectral sequence of modules for the untwisted AHSS.  Specifically, $d_i(ab) = d_i^u(a) b + (-1)^{|a|} a d_i(b)$ where $a$ comes from the untwisted spectral sequence (with untwisted differentials $d_i^u$).
 \end{enumerate} 

All of these properties lift to the twisted differential case with the natural modifications needed in order to make the statement sensible. We will spell this out in detail.

\begin{proposition}
[Properties of the differentials in $\widehat{\rm AHSS}_{\hat \tau}$]
Let $\widehat{\R}$ be a differential ring spectrum refining a ring spectrum 
$\R$. Let $\widehat{\R}_{\hat{\tau}}$ be a differential twist. Then we have the following properties
\begin{enumerate}[label={\bf (\roman*)}]

\vspace{-2mm}
\item {\rm (Linearity)} Each differential 
$d_i$ is a $\widehat{R}^*(\ast)\simeq R(\ast)$-module map.

\vspace{-2mm}
\item {\rm (Normalization)} The twisted differential with a zero twist reduces to the untwisted differential (which may be zero).

\vspace{-2mm}
\item {\rm (Naturality)} If $f:M \to N$ is smooth map between compact manifolds, $f^*$ induces a map of spectral sequences from the 
$\widehat{\rm AHSS}_{\hat \tau}$
 for $(N, \widehat{\R}_{\hat \tau})$ to the 
 $\widehat{\rm AHSS}_{f^*{\hat \tau}}$
  for $(M, f^* \widehat{\R}_{\hat{\tau}})$.  On $E_2$-terms it is induced by $f^*$ in {\v C}ech cohomology, and on 
$E_\infty$-terms it is the associated graded map induced by $f^*$ in the twisted theory $\widehat{\R}_{\hat{\tau}}$.

\vspace{-2mm}
\item {\rm (Module)} The $\widehat{\rm AHSS}_{\hat \tau}$
 for $(N, \widehat{\R}_{\hat{\tau}})$ is a spectral sequence of modules for the untwisted 
 $\widehat{\rm AHSS}$.  Specifically, $d_i(ab) = d_i^u(a) b + (-1)^{|a|} a d_i(b)$ 
 where $a$ comes from the untwisted spectral sequence (with differentials $d_i^u$).
\end{enumerate}
\end{proposition}
\theproof
To prove (i), note that the map $\mu$ in \eqref{mod maprt} gives ${\widehat{R}}(M; \hat{\tau})$ the structure of a module over $\widehat{R}^*(\ast)$ for each $M$ by precomposing with the map induced by the terminal map $M\to \ast$. Clearly this structure is natural with respect to maps of pairs $f:(M,f^*\hat{\tau})\to (N,\hat{\tau})$, with $\tau$ a twist on $M$. Therefore, the maps in the exact couple defining the $\widehat{\rm AHSS}_{\hat \tau}$ commute with the module action. Since the differentials are defined via these maps, the claim follows. 

Property (ii) follows immediately from the observation that the zero twist $\hat{{\tau}}:M\to \ast \to {\rm Tw}_{\widehat{\R}}$ picks out the triple $(\R,c,A)$, where the ring spectrum $\R$ and the algebra $A$ are thought of as modules over themselves. 

Property (iii) also follows directly, as a smooth map $f:M\to N$ induces a morphism of stacks $f^*:{\rm Tw}_{\widehat{\R}}(N)\to  {\rm Tw}_{\widehat{\R}}(M)$. This map sends the sheaf of spectra $\widehat{\R}_{\tau}=(\R_{\tau},t,{\cal L})$ on $N$ to the sheaf of spectra 
$f^*\widehat{\R}_{\tau}=(f^*\R_{\tau},f^*(t),f^*({\cal L}))$ on $M$ via the change of base functor
 $$
 f^*:\sh_{\infty}(M;\Sp)\longrightarrow \sh_{\infty}(N;\Sp)\;.
 $$
 Moreover, the filtration $F_p$ leads to a pullback filtration $f^*F_p$ which is the simplicial diagram formed from the nerve of the cover $\{f^{-1}(U_{\alpha})\}$. This induces the desired morphism of spectral sequences. Considering the relevant diagram that gives the {\v C}ech 
 differential on the $E^1$-page \cite{GS3} and comparing with the corresponding diagram for 
  the cover $\{f^{-1}(U_{\alpha})\}$ one immediately sees that $f$ gives rise to 
  the induced map in {\v C}ech cohomology on the $E_2$-page. 
  
   The proof for property (iv) follows verbatim as in the proof for ring spectrum presented in the differential untwisted case \cite{GS3}, by simply replacing the product map with the module map. 
\endofproof

\section{The AHSS in twisted differential K-theory}
\label{Sec-K}

 By the discussion in the Introduction, and via the general construction in \cite{MSi}, given a cohomology theory $\E$ (represented by an $E_{\infty}$ ring spectrum), the action by automorphisms $\op{GL}_1(\E)$ on $\E$ gives rise to a bundle of spectra over the quotient, which is classified by a map into the delooping $\op{B GL}_1(\E)$. This space classifies the twists of the theory. 
 For K-theory,  by \cite{AGG} ( see also \cite{KS2}), there is essentially a unique equivalence  class of maps of  Picard-$\infty$-groupoids $K(\Z, 2) \to {\rm GL}_1(\mathscr{K})$,  
 which gives the most interesting part of the twists. Delooping the embedding gives a map $\op{B^2U(1)}\into \op{B GL}_1(\mathscr{K})$ and we consider only the twists arising as maps to $\op{B^2U(1)}$, i.e. degree three cohomology classes.

\subsection{Twisting differential $K$-theory by gerbes}
\label{Sec-gerb}

In this section, we give a comprehensive account of the twists in differential $K$-theory 
and discuss the situation from various angles, unifying certain perspectives taken up in 
\cite{BN}\cite{CMW}\cite{KV}\cite{Pa}. 
Our end goal in this section will be to show that differential $K$-theory can indeed be twisted by gerbes, equipped with connection. This is implicit in \cite{BN}, but here we rephrase this result in a way that makes contact with the moduli stack of gerbes $\BB^2U(1)_{\nabla}$ (see \cite{Urs}\cite{FSSt}\cite{SSS3}).\footnote{While 
 this might seem like a further undesirable abstraction, it will be more of a computational and applicable type of abstraction. For example, it will make it
  relatively easy to write down local cocycle data for the twist.}

\medskip
We begin by recounting the topological case. 
Let $\mathscr{K}$ denote the complex K-theory spectrum. The delooping of the units $\op{BGL}_1(\mathscr{K})$ split into three factors. The factor which has attracted the most attention is the Eilenberg-MacLane space $K(\ZZ,3)\into \op{BGL}_1(\mathscr{K})$. Now, recalling that
${\rm Pic}_{\R} \cong \op{BGL}_1(\R)\times \pi_0({\rm Pic}_\R)$,
the inclusion at the identity component of ${\rm Pic}_\R$ then gives rise to a canonical map
\(\label{kz3 tw}
\xymatrix{
\op{B^2U(1)}\simeq K(\ZZ,3) \; \ar@{^(->}[r] & {\rm Pic}_{\mathscr{K}}\;.
}
\)
The left hand side can be identified with the moduli space of bundles with fiber 
$\op{PU}({\cal H})$,  the projective unitary group acting on an 
infinite-dimensional, separable Hilbert space ${\cal H}$.
Indeed, a map $h:M\to \op{B^2U(1)}$ classifies a fibration $\op{BU(1)}\into P\onto M$. The fiber $\op{BU(1)}$ is a classifying space of complex line bundles and has the homotopy type of $\op{PU}({\cal H})$. Alternatively, we can think of $\op{B^2U(1)}$ as the classifying space of \emph{topological} gerbes. Now $\op{PU}({\cal H})$ acts on the space of Fredholm operators ${\cal F}{\rm r}$, which is well known to be a classifying space for $K$-theory (see \cite{At}). Given a map 
$h:M\to \op{B^2U(1)}$, this action leads to an associated bundle of spectra $E_{h}\to M$, which we can regard as a bundle of spectra with fiber $\mathscr{K}$. Taking the homotopy classes of sections of this bundle gives the twisted $K$-theory 
$K(M; h)$. 

\medskip
This is the well-known description for twisted topological $K$-theory established by Atiyah and Segal in \cite{AS1}. In order to make contact with the stacky differential spectra provided by Bunke and Nikolaus, it is useful to take a different perspective. The transition from the above approach to the approach of \cite{BN} is achieved by considering the bundle $E_h\to M$ as a sheaf of spectra on $M$ via its local sections. For an object $(i:N\to M)\in \mf/M$, the value of this sheaf at $N$ is given by the function spectrum $\mathscr{K}_{h}(N):= \Gamma (N, i^*(E_h))$. One can show that the assignment leads to a functor $\mf/M\to \Sp$ and satisfies descent. \footnote{Descent follows from the definition of the coverage via \emph{good} open covers, along with homotopy invariance for $\mathscr{K}$.}
 Since every point on a manifold $M$ admits a geodesically convex neighborhood, this immediately implies that this sheaf of spectra is locally constant and is equivalent to the untwisted K-theory spectrum (regarded as a geometrically discrete sheaf of spectra) locally. As we have seen in Section \ref{Sec-rev} and Section
 \ref{Sec-can}, the actual equivalences which identify the spectrum with untwisted $K$-theory locally can be regarded as a local trivialization of the bundle $E_h\to M$. In this context, descent can be summarized by the diagram in $T({\infty}\mathscr{G}{\rm pd})$
 
\(\label{K bunc}
\hspace{-.3cm}
\xymatrix@C=1.3em{
\hdots \coprod_{\alpha\beta\gamma}  \mathscr{K}\times (U_{\alpha}\cap U_{\beta}\cap U_{\gamma}) \ar[d] \ar@<.1cm>[r]\ar[r]\ar@<-.1cm>[r] &  \coprod_{\alpha\beta}  \mathscr{K}\times (U_{\alpha}\cap U_{\beta}) \ar@<.05cm>[r]\ar@<-.05cm>[r]\ar[d]  & \coprod_{\alpha}  \mathscr{K}\times U_{\alpha}\ar[d] \ar[r] & E_h\ar[d] \ar[r] &  \mathscr{K}/\!/\op{PU}({\cal H})\ar[d]
 \\
\hdots \coprod_{\alpha\beta\gamma}  U_{\alpha}\cap U_{\beta}\cap U_{\gamma}  \ar@<.1cm>[r]\ar[r]\ar@<-.1cm>[r] &  \coprod_{\alpha\beta} U_{\alpha}\cap U_{\beta} \ar@<.05cm>[r]\ar@<-.05cm>[r]  &U_{\alpha} \ar[r] & M\ar[r]^-{h} & K(\ZZ,3)\;,
}
\)
where each square is Cartesian and the top and bottom simplicial diagrams are colimiting. This gives us all the data needed to construct the twisted theory in terms of local data. The simplicial maps in the top diagram are determined by the twist $h:M\to K(\ZZ,3)$ which, more explicitly, is determined by a {\v C}ech 3-cocycle on $M$ via descent.
What we have described so far is the relationship between the bundle of spectra $E_h\to M$ associated to a classifying map $h:M\to K(\ZZ,3)$ and the locally constant sheaf $\mathscr{K}_h$ given by its sheaf of sections. 

\medskip
The next step is to go differential. Crucially, this involves a sort of twisted de Rham theorem which relates twisted $K$-theory to the complex of periodic forms $\big(\Omega^*[u, u^{-1}](M), d_H\big)$ equipped with the differential $d_H=d+H$, where $H$ is a closed 3-form and the differential acts by 
$$
\omega_{2k}u^k+\omega_{2k-2}u^{k-1}\longmapsto 
(d\omega_{2k}+H\wedge \omega_{2k-2})u^k\;.
$$
To understand this interaction, it is important to consider the untwisted case first. 
$K$-theory is an example of a \emph{differentiably simple} spectrum (see \cite{BN}). For such untwisted spectra, there is a rather canonical differential refinement, which in this case is given by taking $A=\pi_*(K)\otimes \RR \cong \RR[u,u^{-1}]$. The usual Chern character map gives an isomorphism
of rings
$$
{\rm ch}:(\mathscr{K}\wedge H\RR)^*(M)
\overset{\cong}{\longrightarrow} H^*\big(\Omega^*[u, u^{-1}](M)\big)\;,
$$
where $\Omega[u, u^{-1}](M)$ is the periodic de Rham complex, equipped with the usual differential. This map refines to a morphism of sheaves of spectra and gives rise to the triple $\widehat{\mathscr{K}}=(\mathscr{K},{\rm ch},\RR[u,u^{-1}])$, which is a differential 
 spectrum representing differential $K$-theory. 

\medskip
Returning to the twisted case, we begin with a closed 3-form $H$, defined on a smooth manifold $M$ and consider the twisted periodic complex $(\Omega^*[u,u^{-1}],d_{H})$ (as a sheaf of DGA's on M).
This complex defines a K-flat, invertible module over the untwisted complex $\Omega^*[u,u^{-1}]$ (see \cite{BN}), with inverse given by $(\Omega^*[u,u^{-1}],d_{-H})$. From the discussion in Section \ref{Sec-flat}, we see that the differential $d_{H}$ should be thought of as a connection on the bundle of spectra $E_h\wedge \Exterior^\bullet(T^*M)\to M$ and the local parallel sections of this connection should give rise to a system of local trivializations of the bundle. We have the following well-known fact which will 
be essential for a detailed description of the trivialization. 
\begin{lemma}
Let $H$ be a closed 3-form on a finite-dimensional smooth manifold $M$ and let $\{U_{\alpha}\}$ be a good open cover of $M$. Let $B_{\alpha}$ be a local potential (i.e. trivialization) for $H$. Then the local sections given by the formal exponentials 
$$
e^{-B_{\alpha}}=1-B_{\alpha}+\tfrac{1}{2!}B_{\alpha}^2 -\hdots 
$$
are twisted closed over $U_{\alpha}$. Moreover, they define isomorphisms of modules via the module action
$$
e^{-B_{\alpha}}\wedge (-):(\Omega_{U_{\alpha}}^*[u,u^{-1}],d)\longrightarrow (\Omega_{U_{\alpha}}^*[u,u^{-1}],d_{H})\;.
$$
\end{lemma}
\theproof
Applying the twisted differential to our section gives
$$
(d+H)(e^{-B_{\alpha}})= d(e^{-B_{\alpha}})+H\wedge e^{-B_{\alpha}}= -dB_{\alpha}e^{B_{\alpha}}+H\wedge e^{-B_{\alpha}}=0\;,
$$
so that $e^{-B_{\alpha}}$ is indeed twisted closed. Now the map $e^{-B_{\alpha}}\wedge (-)$ defines a map of modules since, for each smooth map $f:N\to U_{\alpha}$ 
the pullback form $f^*e^{-B_{\alpha}}$ satisfies
\bea
(d+f^*H)(f^*e^{-B_{\alpha}}\wedge \omega)&=& df^*(e^{-B_{\alpha}})\wedge \omega+f^*(H\wedge e^{-B_{\alpha}})\wedge \omega
\\
&=& f^*(-dB_{\alpha}\wedge e^{-B_{\alpha}})\wedge \omega+f^*(e^{-B_{\alpha}})\wedge d\omega+f^*(H\wedge e^{-B_{\alpha}})\wedge \omega
\\
&=& f^*(e^{-B_{\alpha}})\wedge d\omega\;.
\eea 
Consequently, $f^*(e^{-B_{\alpha}})\wedge (-) $ defines a chain map on $N$ and, as a result, 
$e^{-B_{\alpha}}\wedge (-)$ defines a morphism of sheaves of chain complexes. 
That this map commutes with the module action is obvious from the definition. Moreover, it is easy to check that this map admits an inverse given by $e^{B_{\alpha}}\wedge (-)$. 
\endofproof

\medskip
Thus, we see that the forms $e^{-B_{\alpha}}$ play the role of the $s_{\alpha}$'s in the 
general discussion from Section \ref{Sec-flat}. We depict this in Figure 1.

\begin{figure}[h!]
    \includegraphics[width=0.5\textwidth]{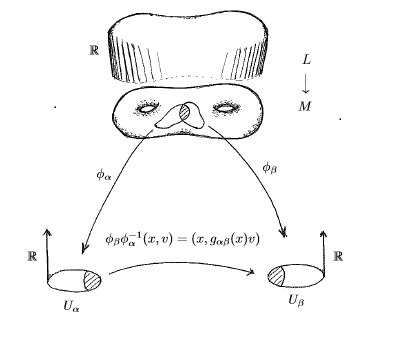} \includegraphics[width=0.5\textwidth]{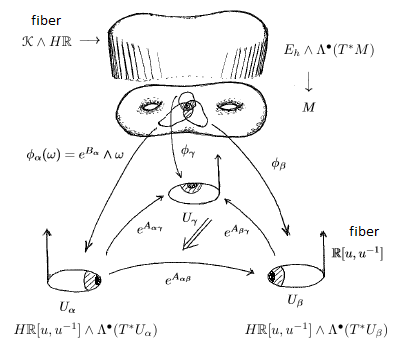}
    
     \caption{Flat line bundle vs. twisted differential cohomology theory.}
     
\end{figure}
\FloatBarrier

\medskip
\noindent Taking $\R$ to be $\mathscr{K}$ in  Theorem \ref{thm-sections}, we get an induced diagram 

\vspace{-4mm}
$$
\xymatrix@R=1.5em{
\vdots & \vdots
\\
\Delta[1]\ar@<.15cm>[u] \ar[u] \ar@<-.15cm>[u]  \ar[r]^-{\eta_{\alpha\beta}} & \prod_{\alpha\beta}B{\rm GL}_1(\mathscr{K}\wedge H\RR)\ar@<.15cm>[u] \ar[u] \ar@<-.15cm>[u] 
\\
\Delta[0] \ar@<.1cm>[u] \ar@<-.1cm>[u]\ar[r]& \prod_{\alpha}B{\rm GL}_1(\mathscr{K}\wedge H\RR) \ar@<.1cm>[u] \ar@<-.1cm>[u]
\;.
}
$$
To see exactly what the maps in this diagram look like, let us choose 1-forms $A_{\alpha\beta}$, defined on intersections of the fixed open cover $\{U_{\alpha}\}$ which satisfy 
$$
dA_{\alpha\beta}=B_{\alpha}-B_{\beta}\;,
$$
and smooth functions satisfying 
$
df_{\alpha\beta\gamma}=A_{\alpha\beta}+A_{\beta\gamma}-A_{\alpha\gamma}
$.
A choice of such data defines a representative for the closed 3-form $H$ in the total complex of the
 {\v C}ech-de Rham double complex. Then we have an induced commutative diagram 
$$
\xymatrix@R=1em{
(\Omega_{U_{\alpha\beta}}^*[u, u^-1],d)\ar[drrrr]^-{\ {\rm exp}(-B_{\alpha})}\ar[dd]_-{{\rm exp}(dA_{\alpha\beta})} &&&&
\\
&&&& (\Omega_{U_{\alpha\beta}}^*[u, u^{-1}],d_H)\;,
\ar[dllll]^-{\ {\rm exp}(B_{\beta})}
\\
(\Omega_{U_{\alpha\beta}}^*[u, u^{-1}],d) &&&&
}
$$
and a similar (3-dimensional) diagram for the triple intersections. Now the equivalence ${\rm ch}:\mathscr{K}\wedge H\RR\to \RR[u,u^{-1}]$ has a pleasant geometric interpretation arising 
from the fact that $K$-theory represents isomorphism classes of vector bundles. Indeed, if we 
take the model for $K$-theory whose zero space is that of Fredholm operators (see \cite{At}), 
then a map
$$
f:M\longrightarrow {\cal F}{\rm r}
$$
defines a virtual vector bundle $(E\ominus F)\to M$ whose fiber at $x\in M$ is given by the virtual difference ${\rm ker}(f(x))-{\rm coker}(f(x))$. If we fix curvature 2-forms ${\cal F}$ and ${\cal G}$ for the bundles $E$ and $F$, respectively, then the Chern character form
$$
{\rm ch}(E\ominus F)={\rm tr}\big({\rm exp}
\big(\tfrac{i}{2\pi}{\cal F}-\tfrac{i}{2\pi}{\cal G}\big)\big)
$$
 constitutes  a geometric representative for 
$
{\rm ch}:\mathscr{K}\wedge H\RR\to \Omega^*[u,u^{-1}]
$,
defined at the level of sheaves of spectra (see \cite{BNV} for a discussion on the cycle map). 
In order to spell out the data guaranteed by Theorem \ref{thm-sections} in the case of K-theory, 
let us consider the task of finding a homotopy commutative diagram 
\(\label{hcom kth}
\xymatrix{
\mathscr{K}\ar[d]^-{\otimes {\cal L}_{\alpha\beta}}\ar[r] & \mathscr{K}\wedge H\RR \ar[rr]^-{\rm ch} \ar[d]^-{\otimes {\cal L}_{\alpha\beta}\wedge {\rm id}} && \Omega_{U_{\alpha\beta}}^*[u,u^{-1}]\ar[rr]^-{e^{-B_{\alpha}}}\ar[d]^-{e^{dA_{\alpha\beta}}} && (\Omega_{U_{\alpha\beta}}^*[u,u^{-1}],d_{H})\ar@{=}[d]
\\
\mathscr{K}\ar[r] &\mathscr{K}\wedge H\RR \ar[rr]^-{\rm ch} &&\Omega_{U_{\alpha\beta}}^*[u,u^{-1}]\ar[rr]^-{e^{-B_{\beta}}}\ar@/^1pc/[ll] && (\Omega_{U_{\alpha\beta}}^*[u,u^{-1}],d_{H})
\;,}
\)
where ${\cal L}_{\alpha\beta}$ is a line bundle defined on intersections, 
giving an action on sections $\Gamma(U_{\a \b};\mathscr{K})$ via 
$$
\big(E_{{\alpha}}\ominus F_{\alpha}\big)\longmapsto \big({\cal L}_{\alpha\beta}\otimes E_{\beta}\ominus {\cal L}_{\alpha\beta}\otimes F_{\beta}\big)\;,
$$
with $E_{\alpha}\ominus F_{\alpha}$ and ${\cal L}_{\alpha\beta}\otimes E_{\beta}\ominus {\cal L}_{\alpha\beta}\otimes F_{\beta}$ being virtual difference bundles on double intersections. 
 After choosing connections on each bundle and, using the fact that the Chern character defines a ring homomorphism, we get the following transformation at the level of forms
\bea
{\rm ch}({\cal L}_{\alpha\beta}\otimes E_{\beta}\ominus{\cal L}_{\alpha\beta}\otimes F_{\beta}) &=& {\rm ch}({\cal L}_{\alpha\beta})\wedge {\rm ch}(E_{\beta})-{\rm ch}({\cal L}_{\alpha\beta})\wedge {\rm ch}( F_{\beta})
\\
&=&{\rm ch}({\cal L}_{\alpha\beta})\wedge \big({\rm ch}(E_{\beta})-{\rm ch}( F_{\beta})\big)\;.
\eea
Here we have
${\rm ch}({\cal L}_{\alpha\beta})={\rm exp}
\big(\tfrac{i}{2\pi}f_{\alpha\beta}\big)
$
and $dA_{\alpha\beta}=\tfrac{i}{2\pi}f_{\alpha\beta}$ a closed 2-form on $U_{\alpha\beta}$, with $f_{\alpha\beta}$ giving the curvature of ${\cal L}_{\alpha\beta}$. Thus, if $H$ represents the curvature of a gerbe then we can choose cocycle data
$$
\widehat{H}=(g_{\alpha\beta\gamma},A_{\alpha\beta},B_{\alpha})\;,
$$
such that $A_{\alpha\beta}$ defines the connection of a line bundle on intersections, whose transition functions $g_{\alpha\beta\gamma}$ on an intersecting patch $U_{\alpha\beta}\cap U_{\gamma}$ satisfy the {\v C}ech cocycle condition on fourfold intersections. In this case, we are able to construct such a homotopy commutative diagram \eqref{hcom kth} and the forms $A_{\alpha\beta}$ and 
$g_{\alpha\beta\gamma}$ define the homotopies and higher homotopies, respectively,
 filling the diagram. We now summarize the above discussion.
 
\begin{proposition}
[Twisted differential K-theory from gerbe data]
Let $\hat{h}:M\to \BB^2U(1)_{\nabla}$ be a gerbe on $M$, with corresponding cocycle data given by the triple $(g_{\alpha\beta\gamma},A_{\alpha\beta},B_{\alpha})$ satisfying the usual gerbe compatibility. Then $\hat{h}$ determines a twisted differential $K$-theory spectrum $\widehat{\mathscr{K}}_{\hat{h}}$,  unique up to a contractible choice.
\end{proposition}

\subsection{The differential K-theory  twisted Chern character}
\label{Sec-ch}

In this section we will consider the effect of rationalization,  
making contact with information encoded in the differential form
part of twisted differential K-theory. To begin with, 
tensoring with the rationals in  the underlying 
topological K-theory gives the  isomorphisms 
$$
(K^0\otimes \QQ)(X) \cong \QQ \times \prod_{k>0} H^{2k}(X; \pi_{2k}(BU) \otimes \QQ)
= H^{\rm ev}(X; \QQ)
\quad \text{and} \quad 
(K^1\otimes \QQ)(X) \cong  H^{\rm odd}(X; \QQ)\;.
$$
The rationalized Chern character is 
$
{\rm ch}_\QQ: K^*(X; \QQ) \to H^*(-; \QQ)
$, so that the Chern character 
can be viewed as the composite map (see \cite{AnH})
$$
\xymatrix{
K^*(X) \ar[r] & K^*(X) \otimes \QQ \cong K^*(X; \QQ) \ar[r]^-{{\rm ch}_\QQ}
& H^*(X; \QQ)
}\;,
$$
where the isomorphism follows from the fact that $\QQ$ is torsion-free.
A key ingredient in the identification of the differentials in the spectral sequence for twisted $K$-theory is  the twisted Chern character map
$$
{\rm ch}_{H}:K(-; {h})\longrightarrow H^{\rm even}(-; H)\;,
$$
which is a natural transformation from $h$-twisted $K$-theory to the $H$-twisted de Rham cohomology of a space, with $H$ a differential form representative for the rationalization of $h$. In order to identify the differentials in the AHSS, we will need a differential refinement of this map. To this end, let us recall 
(see \cite{BSS}\cite{BNV}) that 2-periodic rational cohomology admits a differential refinement via the homotopy pullback
\(\label{Q per diff}
\xymatrix{
\widehat{H\QQ}[u, u^{-1}] \ar[rr] \ar[d]&& 
H(\tau_{\geq 0}\Omega^*\otimes \RR[u,u^{-1}])\ar[d]
\\
H\QQ[u,u^{-1}]\ar[rr]^{i} && H\RR[u,u^{-1}]
\;,
}
\)
where $\tau_{\geq 0}$ denotes the truncation functor which discards all information in negative degrees. This leads to a differential cohomology diamond diagram \cite{SSu} and, in particular, to an exact sequence
$$
\xymatrix{
0\ar[r] & H^{\rm odd}(M;\RR/\QQ)\ar[r] & 
\widehat{H}^{\rm even}(M;\QQ)\ar[r] &
\Omega^{\rm even}_{\rm cl}(M;\RR)\;.
}
$$
If we apply the AHSS to this differential  spectrum, then we get 
the following.
\begin{lemma}
[$\widehat{\rm AHSS}$ for refinement of 2-periodic rational cohomology]
The $E_2$-page for refined 2-periodic rational cohomology  $\widehat{H\QQ}[u, u^{-1}]$  looks as follows
{\small
\(\label{Q per e2ss}
\begin{tikzpicture}
  \matrix (m) [matrix of math nodes,
    nodes in empty cells,nodes={minimum width=5ex,
    minimum height=4ex,outer sep=-5pt}, 
    column sep=1ex,row sep=1ex]{
                &      &     &     & 
                \\         
            1 &     &   &   &
            \\
             0     &   \Omega^{\rm even}_{\rm cl}(M)\oplus \QQ  &  &   &
             \\           
               -1\  \   &  & H^1(M;\RR/\QQ) &  H^2(M;\RR/\QQ)     & 
                \\
              -2\ \  &   &    &  0 & 0 & 0
               \\
               -3\ \ &  & & &  & H^4(M;\RR/\QQ)
               \\
               -4\ \ &   &    &   &  &  & 
               \\
                &      &     &     & \\};
                
  \draw[-stealth] (m-3-2.east) -- node[above]{\small $d_2$} (m-4-4.north west);
   \draw[-stealth] (m-4-3.south east) -- (m-5-5.west);
      \draw[-stealth] (m-5-4.south east) -- (m-6-6.north west);
        \draw[-stealth] (m-4-4.south east) -- (m-5-6.north west);
\draw[thick] (m-1-1.east) -- (m-8-1.east) ;
\end{tikzpicture}
\)
}
\end{lemma}
In particular, the form of the spectral sequence is exactly the same as for $\widehat{K}$, but with $\QQ$'s appearing as the quotients instead of $\ZZ$'s. In this case, all the odd differentials vanish, while the even ones were calculated in \cite{GS3} whose
effect is  shown to extract the class of the singular cocycle \footnote{In \cite{GS3}, the calculation was performed for ordinary differential cohomology, but the same calculation applies to $\QQ$-coefficients.} 
$$
\omega\longmapsto \Big(c\mapsto \oint_{c}\omega_{2k}\mod \QQ \Big)\;,
$$
where $\omega=\omega_0+\omega_2+\cdots \omega_{2k}+\cdots$. There is a differential refinement of the ordinary Chern character map to
$$
\widehat{\rm ch}:\widehat{K}(-)\longrightarrow \widehat{H}^{\rm even}(-;\QQ)\;,
$$
where $\widehat{H}^{\rm even}(-;\QQ)$ is the cohomology theory represented by the differential  spectrum \eqref{Q per diff} (see, \cite{GS3}\cite{Bun} for details). The induced map at the level of spectral sequences for \emph{untwisted} $K$-theory does not yield much information since the only nonvanishing differentials are the even ones, which say something about the values that the Chern character can take on cycles. With rational coefficients, this simply says that the Chern character should take rational values, which is already the case.

\medskip
We now move to the twisted case. Here, we would like to get a differential refinement of the twisted Chern character and see what new conditions arise from the corresponding morphism of AHSS's. To this end, notice that if $\widehat{h}$ is a twist of 
$\widehat{K}$,  then under the map
$$
r:\BB^2U(1)_{\nabla} \cong \BB^2(\RR/\ZZ)_{\nabla} 
\longrightarrow \BB^2(\RR/\QQ)_{\nabla}\;,
$$
induced by the inclusion $r:\ZZ\into \QQ$, we get a twist of $\widehat{H}^{\rm even}(-;\QQ)$. Recall from Theorem \ref{Thm-equiv} that the {\v C}ech cocycle data for $\widehat{h}$ gave the correct gluing conditions for the differential spectrum associated with differential $K$-theory. More precisely, we have a colimiting diagram in $T(\sh_{\infty}(\mf))$, 
$$
\xymatrix{
\hdots   \ar@<.1cm>[r]\ar[r]\ar@<-.1cm>[r]  &  \coprod_{\alpha\beta}  \widehat{\mathscr{K}}\times U_{\alpha\beta}\ar[d]  \ar@<.05cm>[r]\ar@<-.05cm>[r]  &  \coprod_{\alpha}  \widehat{\mathscr{K}}\times U_{\alpha}  \ar[r]\ar[d] & \widehat{\mathscr{K}}_{\widehat{h}}\ar[d]
\\
\hdots   \ar@<.1cm>[r]\ar[r]\ar@<-.1cm>[r]  &  \coprod_{\alpha\beta}  U_{\alpha\beta}  \ar@<.05cm>[r]\ar@<-.05cm>[r]  & \coprod_{\alpha}  U_{\alpha}   \ar[r] & M
\;,
}
$$
where the maps in the diagram are completely determined by the cocycle data for the twist $\widehat{h}:M\to \BB^2U(1)_{\nabla}$. Applying the differential Chern character map locally to each restriction 
$\widehat{\mathscr{K}}\times {U_{\alpha_0\hdots \alpha_n}}$ and using the fact that $\widehat{H}^{\rm even}(-;\QQ)$ is itself twisted by $\widehat{h}$, we get an induced morphism of bundles
\(
\label{chh}
\xymatrix@R=1.2em{
\widehat{\mathscr{K}}_{\widehat{h}}\ar[rr]^-{\widehat{\rm ch}_{\widehat{h}}}\ar[dr] && \widehat{H}\QQ[u,u^{-1}]_{\widehat{h}}
\ar[dl]
\\
& M &
}\;.
\)
\begin{definition}
[Twisted differential Chern character]
We define the twisted differential Chern character 
$\widehat{\rm ch}_{\widehat{h}}$ 
as the map
in diagram\eqref{chh}.
\end{definition}
We would like to get some information about the induced map at the level of spectral sequences. However, we first need to identify what types of differentials arise for 
$\widehat{H}_\QQ^{\rm even}(-; {\hat{h}})$. Since this theory refines 
$H_{\QQ}^{\rm even}(-; h)$, we expect the differentials to be related, and indeed they are. The $E_2$-page for $\widehat{H}_\QQ^{\rm even}(-; {\hat{h}})$
 takes the same form as \eqref{Q per e2ss}, with the only difference occurring at the level of the form part in bidegree $(0,0)$. Here, we have differentials on the $E_{2k}$-page of the form
\(\label{tw oprq}
d_{2k}:(\Omega^{\rm even},d_{H})_{\rm cl}(M)
\longrightarrow
 H^{2k}(M;\RR/\QQ)\;,
\)
where $(\Omega^{\rm even},d_{H})$ is the $H$-twisted de Rham complex and the subscript ${\rm cl}$ indicates that we are taking twisted closed forms in degree zero (i.e. twisted closed even forms). In the twisted case, the maps $d_{2k}$ 
 are modified from the untwisted case. We will explore this in the following sections.

\subsection{Identification of the low degree differentials}
\label{Sec-low}

The formula for the first differential $d_3$ in the twisted AHSS for $K$-theory can be obtained by observing that the differential defines a cohomology operation and that the only cohomology operations (for spaces equipped with maps to $K(\ZZ,3)$) which raise the degree by 3 are given by elements of
$$
H^{k+3}(K(\ZZ,k)\times K(\ZZ,3);\ZZ)\cong H^{p+3}(K(\ZZ,k);\ZZ)\oplus H^{p+3}(K(\ZZ,3);\ZZ)\oplus \ZZ\;,
$$
where the third summand is generated by the product of the generators for $H^k(K(\ZZ,k);\ZZ)$ and $H^3(K(\ZZ,3);\ZZ)$. This implies the formula 
$$
\xymatrix{
d_3=Sq^3_\Z + (-) \cup \lambda h: H^p(X; \Z)=E_{3}^{p,q} 
\; \ar[r] &
\;  E_{3}^{p+3, q-2}=H^{p+3}(X; \Z)\;,
 }
$$
where $\lambda$ is an integer which has yet to be determined. To compute this integer, it is sufficient to consider the spectral sequence on the sphere $S^3$, where one computes $\lambda=-1$ (see \cite{AS2}). 

\medskip
We would like to apply similar reasoning in the differential setting. To do this, however, we should take some care. The AHSS for differential cohomology theory has the additional datum of the level at which we are computing the theory. For example, in the case of $K$-theory, we have a different spectral sequence for each of the differential  spectra
$$
\widehat{\mathscr{K}}_{(n)}:=\big(\Sigma^nK,{\rm ch},\RR[u , u^{-1}][n]\big)\;,
$$
where we have shifted both the $K$-theory spectrum and the complex $\RR[u, u^{-1}]$ by $n$. The degree zero cohomology of this differential  spectrum computes the degree $n$ differential $K$-theory of the underlying manifold. This is not really too surprising since the spectral sequence is based off of the Mayer-Vietoris sequence (as is the AHSS for ordinary $K$-theory). In \cite{BNV}, the Mayer-Vietoris sequence is considered and this sequence depends on the degree of the cohomology theory in question. Therefore, we will have to differentially refine one degree at a time. 
For illustration, we will consider the case $n=0$ as this will compute the differential refinement of $K^0$. Similar effects hold for $K^1$, which we will consider in 
Section \ref{Sec-ex}. 

\medskip
Suppose we are given a twist $\widehat{h}:M\to \BB^2U(1)_{\nabla}$
for the theory $\widehat{K}^0$. Setting $\widehat{\R}=\widehat{\mathscr{K}}$
in Theorem \ref{Thm-ourAHSS}, we observe that the $E^{0,0}_2$-entry in the spectral sequence 
can be calculated as
$$
H^0(M;\widehat{K}(-;\hat{h}))=\ker\Big\{\xymatrix{\prod_{\alpha}\ZZ\times \Omega_{>0}^{\rm even}(U_{\alpha})\ar[rr]^-{e^{dA_{\alpha\beta}} -{\rm id}} && \prod_{\alpha\beta}\ZZ\times \Omega_{>0}^{\rm even}(U_{\alpha\beta})
}\Big\}\;,
$$
where the subscript $>0$ indicates that we are taking forms with vanishing degree-zero term. 
Note that this group is isomorphic to the subgroup of twisted closed forms, whose elements are of the type 
$\omega=\omega_0+\omega_2+\hdots \;,$ with $\omega_0$ an integer. If the curvature of $\hat{h}$ is nontrivial in de Rham cohomology, then $\omega_0$ must vanish (otherwise we could trivialize $H$ by the form $\omega_2/\omega_0$). In this case, the $E^{0,0}_2$-entry is precisely the twisted closed forms and the full $E_2$-page looks as follows
{\small
\(
\label{E2-Kth}
\begin{tikzpicture}
  \matrix (m) [matrix of math nodes,
    nodes in empty cells,nodes={minimum width=5ex,
    minimum height=4ex,outer sep=-5pt}, 
    column sep=1ex,row sep=1ex]{
                &      &     &     & 
                \\         
            1 &     &   &   &
            \\
             0     &   (\Omega^*[u, u^{-1}],d_H)(M)_{\rm cl} &  &   &
             \\           
               -1\  \   &  & H^1(M;U(1)) &  H^2(M;U(1))     & 
                \\
              -2\ \  &   &    &  0 & 0 & 0
               \\
               -3\ \ &  & & &  & H^4(M;U(1))
               \\
               -4\ \ &   &    &   &  &  & 
               \\
                &      &     &     & \\};
                
  \draw[-stealth] (m-3-2.east) -- node[above]{\small $d_2$} (m-4-4.north west);
   \draw[-stealth] (m-4-3.south east) -- (m-5-5.west);
      \draw[-stealth] (m-5-4.south east) -- (m-6-6.north west);
        \draw[-stealth] (m-4-4.south east) -- (m-5-6.north west);
\draw[thick] (m-1-1.east) -- (m-8-1.east) ;
\end{tikzpicture}
\)
}
while the $E_3$-page looks as
{\small
$$
\begin{tikzpicture}
  \matrix (m) [matrix of math nodes,
    nodes in empty cells,nodes={minimum width=5ex,
    minimum height=5ex,outer sep=-5pt}, 
    column sep=1ex,row sep=1ex]{
                &      &     &     & 
                \\         
            1 &     &   &   &
            \\
             0     &   \ker(d_2)  &  &   &
                \\
              -1\ \  &   &    &  & H^2(M;U(1))  & H^3(M;U(1)) 
               \\
               -2\ \ &  &  &  & 0 
               \\
               -3\ \ &   &    &   &  &  & & H^5(M;U(1))
               \\
                 &      &     &     & \\};
                 
             \draw[-stealth] (m-3-2.south east) --node[pos=.5,below]{\small $d_3$} (m-5-5.north west);      
   \draw[-stealth] (m-4-5.south east) --node[pos=.5,below]{\small $d_3$} (m-6-8.north west);
\draw[thick] (m-1-1.east) -- (m-7-1.east) ;
\end{tikzpicture}
$$
}
This pattern continues and we are led to the following. 

\begin{lemma}
[Types of differentials in $\widehat{\rm AHSS}_{\hat h}$ for 
$\widehat{K}^0(-; \hat{h})$]
\label{Lem-2d}
There are two types of differentials:
\begin{enumerate}[label=(\roman*)]

\vspace{-2mm}
\item \emph{(Obstructions associated to curvature forms):} Those of the form
\(\label{type 1 diff}
d^{\rm curv}_{2k}:S\subset 
(\Omega^*[u, u^{-1}],d_{H})(M)_{\rm cl} \longrightarrow
 H^{2k}(M;U(1))\;,
\)
where $S=\bigcap {\rm ker}(d_{2j})$ for $k<j$, which emanate from the entry $E^{0,0}_{2k}$. 

\vspace{-2mm}
\item \emph{(Obstructions associated to flat classes):} 
We also have differentials occurring on the odd pages of another type, namely
\(\label{type 2 diff}
d^{\rm fl}_{2k+1}:H^{p}(M;U(1)) \longrightarrow H^{p+2k+1}(M;U(1))\;,
\)
emanating from the entries $E^{p,-2q+1}_{2k+1}$, with $q\geq 0$. 
\end{enumerate}
\end{lemma}

As the notation suggests, the differentials of the first type gives obstructions for the curvature forms to lift to differential $K$-theory, while the differentials of the second type give rise to obstructions to lifting to flat classes (see e.g. \cite{Lo} \cite{Ho2} for 
a description of such classes).

\medskip
In \cite{GS2}, we classified all differential cohomology operations by computing the connected components of the mapping space $\map(\BB^nU(1)_{\nabla},\BB^mU(1)_{\nabla})$. We showed that all such operations either arise as an integral multiple of a power of the Dixmier-Douady class (in the sense of \cite{cup} \cite{E8})
$$
{\rm DD}:={\rm id}:\BB^nU(1)_{\nabla}\longrightarrow
 \BB^nU(1)_{\nabla}\;,
$$
or arise from an operation of the form $\phi:H^n(-;\ZZ)\to H^{m}(-;U(1))$, via the natural map $I:\widehat{H}^n(-;\ZZ)\to H^n(-;\ZZ)$,
 and the inclusion of flat classes 
 $j: {H}^{m}(-;U(1))\to \widehat{H}^{m+1}(-;\ZZ)$. 
In particular, we can get transformations of the second type by lifting integral cohomology operations through the Bockstein for the exponential sequence $\beta:H^{m}(-;U(1))\to H^{m+1}(-;\ZZ)$. In fact, such lifts are exactly the refinements of the torsion operations in integral cohomology. This allows us to compute the differentials of the second type \eqref{type 2 diff}. 

\begin{proposition}
[Flat first differential for twisted differential K-theory]
\label{Prop-flat1st}
Let $\widehat{h}:M\to \BB^2U(1)_{\nabla}$ be a twist for differential $K$-theory. Then we have
$$
d^{\rm fl}_{3}=\widehat{Sq}^3+ \widehat{h}\cup_{{}_{\rm DB}} (-) \;,
$$
where $\widehat{Sq}^3$ is the operation 
$
\Gamma_2Sq^2\rho_2\beta:H^{p-1}(-;U(1)) \to H^{p+2}(-;U(1))
$ 
and $\cup_{{}_{\rm DB}}$ is the Deligne-Beilinson cup product. 
\footnote{Note that we are considering cohomology with $U(1)$ coefficients 
as sitting inside differential cohomology by the inclusion $j$.}
\end{proposition}
\theproof
By the characterization in \cite{GS2}, the only operation which can raise the degree by 3 comes  from integral multiples of powers of the Dixmier-Douady class or torsion operations in integral cohomology. Moreover, in \cite{GS3}, we showed that the differentials for $\widehat{\rm AHSS}$ must refine the differentials in the underlying AHSS, in the sense that we have commutative diagrams, one for each $p$  
$$
\xymatrix@=1.5em{
H^{p-1}(M;U(1))\ar[rr]^-{\beta}\ar[d]^-{d^{\rm fl}_3} &&
 H^{p}(M;\ZZ)\ar[d]^-{d_3}
\\
H^{p+2}(M;U(1))\ar[rr]^-{\beta} &&
 H^{p+3}(M;\ZZ)\;.
}
$$
Since the differentials in the underlying theory are given by the cup product with the integral class $h$ and $\beta(\widehat{h}\cup_{\rm DB} \hat{\alpha})=h\cup \beta(\hat{\alpha})$, for a flat class $\hat{\alpha}$. In \cite{GS3}, we showed that $\widehat{Sq}^3$ refines the differential $d_3$ in untwisted $K$-theory. Thus, the only possibility is
$d_3=\widehat{Sq}^3+\widehat{h} \cup_{{}_{\rm DB}}(-)$.
\endofproof

We now turn our attention to the differentials $d_p^{\rm curv}$, the first of the two types in Lemma \ref{Lem-2d}. Let us consider the sheaf of spectra $\tau_{\geq p}\QQ[u,u^{-1}]\otimes \Omega^*$, where $\tau_{\geq p}\QQ[u,u^{-1}]$ is the truncated complex, with degrees $\leq p$ replaced by zeros. We have natural maps
$$
\xymatrix{
i:\tau_{\geq p}\QQ[u,u^{-1}]\;\ar@{^{(}->}[r] &
 \tau_{\geq p}\QQ[u,u^{-1}]\otimes \Omega^*\;, 
&
 j:\Omega^*\;\ar@{^{(}->}[r] & 
 \tau_{\geq p}\QQ[u,u^{-1}]\otimes \Omega^*\;.
 }$$
We denote the pullback of $i$ by $j$ in sheaves of spectra by $F_p\widehat{H}\QQ[u,u^{-1}]$ and we have a canonical map
$$
F_p\widehat{H}\QQ[u,u^{-1}]\longrightarrow \widehat{H}\QQ[u,u^{-1}]\;,
$$
which is induced by the inclusion $\tau_{\geq p}\QQ[u,u^{-1}]\into \QQ[u,u^{-1}]$. Since a choice of twist $\hat{h}:M\to \BB^2(\RR/\QQ)_{\nabla}$ raises the degree, we can consider the complex $\tau_{\geq p}\QQ[u,u^{-1}]\otimes \Omega^*$ with differential $d+u H\wedge$, 
\footnote{Note that we have chosen to work with chain complexes rather than in 
cochain complexes, so that differential forms are weighted with negative degrees.}
and the bundle of $H\big(\tau_{\geq p}\QQ[u,u^{-1}]\big)$-spectra ${\cal H}\to M$, induced by the topological twist $h:M\to K(\QQ,3)$. The local exponential map still gives rise to a twisted Chern character and thus we have an $\hat{h}$-twisted theory 
$F_p\widehat{H}\QQ[u,u^{-1}]_{\hat{h}}$, which comes equipped with a canonical map
$$
F_p\widehat{H}\QQ[u,u^{-1}]_{\hat{h}}\;\longrightarrow
 \widehat{H}\QQ[u,u^{-1}]_{\hat{h}}\;,
 $$
induced by the inclusion $i$. This map gives rise to a morphism of spectral sequences which on the $E_2$-page determines a commutative diagram 
\(\label{filt pdiag q}
\xymatrix{
F_p(\Omega^*[u,u^{-1}],d_{H})(M)_{\rm cl}\;
\ar@{^{(}->}[rr]\ar[d]^-{d^{\prime}_{p}} && (\Omega^*[u,u^{-1}],d_{H})(M)_{\rm cl}\ar[d]^-{d_{p}}
\\
{\rm ker}(d^{\prime}_{p-1})\ar[rr] && \ker(d_{p-1})/{\rm im}(d_{p-1})
\;.
}
\)
The spectra $F_p\widehat{H}\QQ[u,u^{-1}]_{\hat{h}}$ give rise to a filtration of the spectrum $\widehat{H}\QQ[u,u^{-1}]_{\hat{h}}$, and we also have an exact sequences of sheaves of spectra
\(\label{exfilp}
\xymatrix{
F_{p+1}\widehat{H}\QQ[u,u^{-1}]_{\hat{h}}\ar[r]&
 F_p\widehat{H}\QQ[u,u^{-1}]_{\hat{h}}\ar[r] &
  F_p\widehat{H}\QQ[u,u^{-1}]_{\hat{h}}/F_{p+1}\widehat{H}\QQ[u,u^{-1}]_{\hat{h}}
}
\)
induced by the inclusion $\tau_{\geq p+1}\QQ[u,u^{-1}]\into \tau_{\geq p}\QQ[u,u^{-1}]$. The rightmost map in \eqref{exfilp} gives rise to a morphism of spectral sequences and, consequently, to  a commutative diagram 
\footnote{This requires a quick calculation, but essentially the point is that the subcomplexes $F_p(\Omega^*[u,u^{-1}],d_{H})(M)_{\rm cl}$ contain those twisted closed forms of the form $\omega=0+0+\hdots \omega_p+\omega_{p+2}+\hdots$. Thus, at the level of forms, the quotient $F_p(\Omega^*[u,u^{-1}],d_{H})(M)_{\rm cl}/F_{p+1}(\Omega^*[u,u^{-1}],d_{H})(M)_{\rm cl}\cong \Omega^p_{\rm cl}(M)$.}
$$
\xymatrix{
F_p(\Omega^*[u,u^{-1}],d_{H})(M)_{\rm cl}\ar[rr]\ar[d]^-{d^{\prime}_{p}} && \Omega^p_{\rm cl}(M)\ar[d]^-{d^{\prime\prime}_{p}}
\\
{\rm ker}(d^{\prime}_{p-1}) \; \ar@{^{(}->}[rr] && H^{p}(M;\RR/\QQ)
\;.
}
$$
The differential $d^{\prime\prime}_{p}$ was identified in \cite{GS3} as the map which takes the class of $\omega_p$ in singular cohomology and mods by $\QQ$. Thus we identify the differential $d^{\prime}_p$ as the map
$$d^{\prime}_p(0+0+\hdots+\omega_p+\omega_{p+2}+\hdots)=[\omega_p]\mod
 \QQ\;.$$
 Returning to the diagram \eqref{filt pdiag q}, we see that when restricted to the forms in the $p$-th level of the filtration $F_p(\Omega^*[u,u^{-1}],d_{H})(M)_{\rm cl}$, $d_{p}$ must give the class of the leading term, modulo $\QQ$. This gives the following characterization of the differentials $d^{\rm curv}_p$, at least when restricted to forms in the $p$-th level of the filtration on $(\Omega^*[u,u^{-1}],d_{H})$.

\begin{proposition}
[Characterization of the differentials $d^{\rm curv}_p$]
Let $p$ be an even integer. For the refinement of periodic rational cohomology, the differential $d_p^{\rm curv}$ take the form
$$
\xymatrix{
d_p: (\Omega^*[u,u^{-1}],d_{H})(M)^{\QQ}_{\rm cl} \; \ar[r] &
  H^{p}(M;\RR/\QQ)
}\;,
$$
where $(\Omega^*[u,u^{-1}],d_{H})(M)^{\QQ}_{\rm cl}$ is the subgroup of twisted closed forms with degree zero term $\omega_0$ given by a constant function taking values in $\QQ$. Moreover, the differential $d_p$ maps a twisted closed form of the type
$\omega=0+0+\hdots+\omega_p+\omega_{p+2}+\hdots$
to the class of the leading term $\omega_p$, modulo $\QQ$, i.e.
$$d_p(0+0+\hdots +\omega_p+\omega_{p+2}+\hdots)=[\omega_p]\mod \QQ\;.$$
For $p=2$ there are two possibilities. Either $\omega_0=0$, or the curvature of the twist $\hat{h}$ is globally exact. 
\begin{enumerate}[label={\bf (\roman*)}]

\vspace{-2mm}
\item  If $\omega_0=0$ then 
$$
d_2(0+\omega_2+\omega_4 \hdots )=[\omega_2]\mod \QQ\;.
$$

\vspace{-2mm}
\item If $H=0$, the twist $\hat{h}$ is trivial, so that we have an 
 isomorphism of spectral sequences between $\widehat{\rm AHSS}_{\hat h}$
 and $\widehat{\rm AHSS}$. This, in turn, gives the 
 commutative diagram
$$
\xymatrix{
(\Omega^*[u,u^{-1}],d_{H})(M)^{\QQ}_{\rm cl}\ar[rr]^-{e^{-B}\wedge (-)}
\ar[d]^-{d^{\prime}_{2}} && \Omega^*[u,u^{-1}](M)_{\rm cl}\ar[d]^-{d_{2}}
\\
H^{2}(M;\RR/\QQ) \ar@{=}[rr] && H^{2}(M;\RR/\QQ)
\;,
}
$$
so that 
$
d_2(\omega)=[\omega_2-B\omega_0]\mod \QQ
$.
\end{enumerate}
\end{proposition}
\begin{remark}
Identifying the differentials $d^{\rm curv}_p$ in full generality seems to be complicated. The structure appears to be related to lifting the terms of the locally defined closed form $e^{-B_{\alpha}}\wedge \omega$ to a flat bundle, using the gerbe data for the twist $\hat{h}$, but we are unsure of how to present the differential in such a way that the lifting condition is useful in practice (i.e. a checkable condition on the twisted forms so that they lift to curvatures of twisted differential $K$-theory). Moreover, we have only discussed the rational case here and the situation for twisted differential $K$-theory seems to be further complicated by the presence of torsion. Since we believe such a presentation exists and would be extremely useful in applications, we save this for future investigation.
\end{remark}

\subsection{The higher flat differentials}
\label{Sec-high}

In \cite{AS2}, Atiyah and Segal considered the spectral sequence for the twisted de Rham complex $(\Omega^*[u,u^{-1}](M),d_{H})$ on a smooth manifold $M$. This complex admits a filtration \footnote{Note that this is the same filtration we considered in the previous section.} 
with the $p$-th level $F_p(\Omega^*[u,u^{-1}](M),d_{H})$ given by the complex which 
in degree $k$ has elements of the form
$$
\omega=0+0+\hdots \omega_p+\omega_{p+2}+\hdots\;.
$$
The spectral sequence with respect to this filtration has $E_2$-page given
by  the ordinary de Rham cohomology groups
$$
E^{p,q}_2=H^{p+q}_{\rm dR}(M)\;.
$$
Atiyah and Segal showed that the differentials $d_{2k+1}$ in the spectral sequence for the periodic complex of forms are given by the Massey product operation $\langle H,H,\hdots,H,x\rangle$, where $H$ is repeated $k$-times. 

\medskip
We highlight that the Massey products are related to solving the twisted closed equation
 $dx=H\wedge x$ in the following way. Let $x_2$ be a degree two de Rham class and 
 consider the Massey triple product $\langle H,H,x_2\rangle$. A representative for the product 
 is constructed by trivializing each of the wedge products $H\wedge H=0$ and $H\wedge x_2$. Clearly, we can choose the zero form as a trivialization of the first product. Assuming that the second product vanishes in cohomology, then we can find $x_4$ such that $dx_4=-H\wedge x_2$. In this case, we have the representative of the Massey triple product
$$
\langle H,H,x_2\rangle=H\wedge x_4\;,
$$
which is a class in cohomology of degree seven. In general, there is some ambiguity in defining the Massey product coming from the choice of $x_4$. This ambiguity is measured precisely by classes 
in $[H]\wedge H^4_{\rm dR}(M)$. Thus, the product is well-defined in the quotient
$$
d_5[x_2]=- \langle H,H,x_2\rangle\in H^7_{\rm dR}(M)/[H]\wedge H^4_{\rm dR}(M)\;,
$$
and this is precisely the quotient one finds on the $E_5$-page of the spectral sequence. Thus, if $d_5$ vanishes on $[x_2]$, then there is $x_6$ such that $H\wedge x_4=dx_6$ and we can form the Massey product $H\wedge x_6$. Again the ambiguity is measured by the subgroup 
$[H]\wedge H^6_{\rm dR}(M)$. If one iterated the process, we find that we are actually asking for a solution to the twisted equation $dx_{2k+2}=-H\wedge x_{2k}$. 

\medskip
We now want to find the relationship between the Massey products on the higher differentials in the spectral sequence for twisted differential $K$-theory. In \cite{GS1} we established the basic theory of 
differential Massey products. Indeed, using these we have the following.

\begin{proposition}
[The higher differentials in $\widehat{\rm AHSS}_{\hat h}$ of $\widehat{K}^0(-; \hat{h})$]
Let $\widehat{h}:M\to \BB^2U(1)_{\nabla}$ be a twist for differential $K$-theory, regarded as a twist for periodic rational cohomology via the differential Chern 
character map (see Section \ref{Sec-ch}). Then the differentials $d_{2k+1}$ can be identified with the differential Massey product operation 
$$
d_{2k+1}= -  \langle \underbrace{ \hat{h},\hat{h}, \hdots, \hat{h}}_{k-{\rm times}} , -\rangle\;.
$$
\end{proposition}
\theproof
Consider the shifted twisted theory $H\big(\RR[u,u^{-1}]\big)[1]_h$ (i.e. the shift of periodic real cohomology). By restriction, $\hat{h}$ twists the underlying flat theory for 
$\widehat{H}\big(\QQ[u,u^{-1}]\big)$, which is equivalent to the shifted spectrum $H\big(\RR/\QQ[u,u^{-1}]\big)[1]$. The equivalence allows us to twist periodic $\RR/\QQ$-cohomology by $\hat{h}$. The mod $\QQ$ map thus gives a well-defined morphism
$$
H\big(\RR[u,u^{-1}]\big)[1]_h\longrightarrow 
H\big(\RR/\QQ[u,u^{-1}]\big)[1]_{\hat{h}}\;,
$$
which induces a morphism of corresponding spectral sequences. At the level of spectral sequences,
 the shift of $H\big(\RR[u,u^{-1}]\big)_h$ manifests itself by simply shifting the entries on each 
 page. By Atiyah and Segal \cite[Prop. 6.1]{AS2}, the differentials for this spectral sequence 
 are the Massey products 
$$
d_{2k+1}[x]=-\langle \underbrace{h,h,\hdots}_{k-{\rm times}},[x] \rangle\;.
$$ 
Since we have a morphism of spectral sequences, it follows that the differentials 
$d_{2k+1}^{\prime}$ on the $(2k+1)$-page of the spectra sequence for 
$H\big(\RR/\QQ[u,u^{-1}]\big)[1]_{\hat{h}}$ refines the Massey product 
operation in that we have a commutative diagram
$$
\xymatrix{
H^p(M;\RR)\ar[rr]^-{\mod \QQ}\ar[d]^{-\langle h,h,\hdots,-\rangle}  && H^p(M;\RR/\QQ)\ar[d]^-{d_{2k+1}^{\prime}}
\\
H^{p+2k+1}(M;\RR)\ar[rr]^-{\mod \QQ} &&  H^{p+2k+1}(M;\RR/\QQ)
\;.
}
$$
Using the formula for the Deligne-Beilinson cup product on flat classes (see \cite{cup}
or expression \eqref{DB-ab}), 
it is easy to see that the differential Massey product operation 
$-\langle \hat{h},\hat{h},\hdots,-\rangle$ indeed refines the underlying 
Massey product. 

We claim that this is the only possible operation.
Indeed, let 
\(
\label{map-phi}
\phi:H^p(M;\RR/\QQ)\longrightarrow
 H^{p+2k+1}(M;\RR/\QQ)
\)
 be a natural operation on manifolds equipped with maps $M\to \BB^2\RR/\QQ_{\nabla}$, refining the Massey product operation in cohomology with real coefficients. Then, from the diagram,
  the difference $\phi+\langle \hat{h},\hat{h},\hdots,-\rangle$, when restricted to $H^p(M;\RR)$, must vanish. But the mod $\QQ$ map is surjective since the Bockstein
$$
\xymatrix{
H^{p}(M;\RR/\QQ)\ar[r]^{\beta_{\QQ}} &
H^{p+1}(M;\QQ)\; \ar@{^{(}->}[r] & H^{p+1}(M;\RR)
}
$$
vanishes. Thus, $\phi =-\langle \hat{h},\hat{h},\hdots,-\rangle$
on $H^p(M;\RR/\QQ)$.
\endofproof

\medskip
We will now show that the corresponding Massey products in 
the topological and the differential cases are compatible
under the composite map 
\(
\label{map-varphi}
\xymatrix{
\varphi: H^{*-1}_{\rm dR}(-; \RR) 
\ar[rr]^-{\mod \QQ} &&
H^{*-1}(-; \RR/\QQ) \ar[r]^-{j}
& 
\widehat{H}^*(-; \QQ)\;.
}
\) 

\begin{lemma}
[Correspondence of Massey products]
\label{difmas hw}
Let $\hat{h}:M\to \BB^nU(1)_{\nabla}$ be a differential cohomology class
with curvature $H$ and let $\omega$ be a degree $m$ 
closed differential form on $M$. When the Massey product  is defined, then, 
for each representative 
$$
x\in \langle  \underbrace{H,H, \hdots, H}_{k-{\rm times}}, \omega \rangle
\subset H^{m+k(n-1) +1}_{\rm dR}(M)
\;,
$$
the map \eqref{map-varphi} gives a corresponding representative 
$$
\varphi(x)\in \langle  \underbrace{\hat{h},\hat{h}, \hdots, \hat{h}}_{k-{\rm times}}, 
\varphi(\omega) \rangle \subset \widehat{H}^{m+k(n-1) +2}(M; \QQ)
\;.
$$
\end{lemma}
\theproof
Recall that we have a formula for the Deligne-Beilinson cup product at the 
level of cochains with values in the Deligne complex with $\QQ$-coefficients. 
At the level of local sections $\alpha=(\alpha_0,\alpha_1,\hdots, \alpha_n)$ and $\beta=(\beta_0,\beta_1,\hdots ,\beta_m)$, this product takes the form
(see \cite{cup}\cite{GS3})
\(
\label{DB-ab}
\alpha\cup_{\rm DB} \beta=\left\{\begin{array}{cc}
\alpha_0\beta_\ell, & \ell \neq m
\\
\alpha_p\wedge d\beta_m ,& p\neq 0.
\\
0, & \text{otherwise}. 
\end{array}
\right.
\)
If $\omega$ is a globally defined closed form, then $\varphi([\omega])$ is the class of 
the cocycle in Deligne cohomology $\varphi(\omega)=(0,0,\hdots,\omega)$.
 On $\hat{h}$ the above formula gives
$$
\varphi(\omega)\cup_{\rm DB} \widehat{h}= \varphi(\omega \wedge H)\;.
$$
Since the DB cup product defines a graded commutative product at the level of cohomology, the above cocycle represents the class $\varphi([\omega])\cup_{\rm DB} {\hat h}$. Since $\varphi$ arose from 
a morphism of DGA's, it follows immediately that $\varphi$ must send a defining system for the Massey product to a defining system for the differential Massey product.
\endofproof

\subsection{Examples} 
\label{Sec-ex}

The basic properties listed above are enough to get started with computations. We illustrate with the example of the 3-sphere $S^3$, which generalizes to differential twisted K-theory the corresponding twisted K-theory calculation in 
\cite{MMS}(after M. Hopkins)\cite{BCMMS}\cite{Br}\cite{Do}\cite{FHT}\cite{MR}. 
The twisted differential K-theory of the Lie group $\op{SU}(2)$ (and also $\op{SU}(3)$)
was studied in \cite{CMW} using index and group theoretic methods. 

\medskip
Let $h:S^3\to K(\ZZ,3)$ be a map representing an integer (which we also denote by $h$) in $H^3(S^3;\ZZ)\cong \ZZ$. We denote the corresponding twisted $K$-theory 
 on $S^3$ by ${K}^*(S^3; {h})$. 
Now the map $h$ can be refined to a gerbe with connection $\widehat{h}:S^3\to \BB^2U(1)_{\nabla}$, with curvature form $H$. Let $(\Omega^*[u,u^{-1}](S^3),d_{H})$ denote the periodic, $H$-twisted de Rham complex on $S^3$, with differential $d+H\wedge (-)$.  
Thus, the triple $\widehat{\mathscr{K}}_{{\hat h}}:=(\mathscr{K}_{h}, {\rm ch}, (\Omega^*[u,u^{-1}],d_{H}))$ gives the data of a differential refinement of $h$-twisted $K$-theory 
spectrum. We denote the underlying theory by $\widehat{K}^*(S^3; \hat{h})$.

\begin{proposition} [Twisted differential $K^0$  of the 3-sphere]
We have
$$
\widehat{K}^0(S^3; {\hat h})\cong \Omega_{\rm cl}^2(S^3)\;.
$$
\end{proposition} 

\theproof
Consider the cover $\{U,V\}$, given by removing the north and south poles of $S^3$, 
respectively. Then the Mayer-Vietoris sequence from Prop. \ref{prop twdiff} (iii)
gives an exact sequence
\(\label{red seqtk}
\xymatrix@=1.5em{
K_{\RR/\ZZ}^{-2}(S^2; \hat{h})\ar[r] & 
\widehat{K}^0(S^3; \hat{h})\ar[r] & 
\widehat{K}^0(U; \hat{h})\oplus \widehat{K}^0(V; \hat{h})\ar[r]^-f &
 \widehat{K}^0(U\cap V; \hat{h})\ar[r] & K^{1}(S^3; h)
}\;.
\)
The twist trivializes on each local patch $U$ and $V$. Moreover, since the restriction of $h$ to $U\cap V\simeq S^2$ trivializes by dimension, the twisted differential $K$-theory of $U\cap V$ also reduces to the untwisted differential $K$-theory. Our job is now reduced to calculating the untwisted groups and understanding the maps in \eqref{red seqtk}. We have
\begin{eqnarray} 
\widehat{K}^0(U) &\cong& \ZZ\oplus \Omega^{\rm even}_{\rm cl}(U) \label{u diffk}
\\
\widehat{K}^0(V) &\cong& \ZZ\oplus \Omega^{\rm even}_{\rm cl}(V) \label{v diffk}
\\
\widehat{K}^0(U\cap V)&\cong& \ZZ\oplus \ZZ\oplus d\Omega^{\rm odd}(U\cap V)\;, \label{int diffk}
\end{eqnarray}
where the first factor of $\ZZ$ in \eqref{int diffk} is generated by $1$ and the second factor is generated by the canonical line bundle $L$ on $S^2$. Now the maps in the sequence depend on the twist because the identification of the twisted theory with the untwisted theory depends on a choice of trivialization, and that trivialization is not canonical (i.e. it depends on the $U$ and $V$). The transition maps between the trivializations give the correction from the untwisted case.

Notice that the factor $\ZZ\oplus \ZZ$ in the last group \eqref{int diffk} correspond to the $\ZZ$-linear span of each of the basis elements $\{1,L\}$ in $K(S^2)\cong \ZZ[L]/(1-L)^2$. The copies of $\ZZ$ in the previous two equations \eqref{u diffk} and \eqref{v diffk} correspond to the $\ZZ$-linear span of the trivial bundle. Topologically, the twist $h$ restricts on the intersection $U\cap V\simeq S^2$ to a cocycle representing $L$. Geometrically, the group structure on $H^3(S^3;\ZZ)\cong \ZZ$ corresponds to the group generated by $L$ under tensor product of line bundles. The action of the transition functions for the local trivializations corresponds to multiplication by $L^{\otimes h}$, which, when written in the basis $\{1,L\}$, is $(h-1)1-hL$. Thus, if we restrict the map $f$ in \eqref{red seqtk} to the copies of $\ZZ$, we get the map
$$
\left(\begin{array}{cc}
1 & h-1
\\
0 & -h
\end{array}
\right):\ZZ\oplus \ZZ\longrightarrow \ZZ\oplus \ZZ\;.
$$
This map is coming from the topological part of twisted $K$-theory. 

We still need to identify the geometric part. To this end, recall that we also have local trivializations for the twisted complex of forms, coming from the formal exponentials $e^{B_{U}}$ and $e^{B_{V}}$, where $B_{U}$ and $B_{V}$ are local potentials (i.e. trivialization) for $H$. Using the transition data $e^{dA}$ on intersections (i.e. $A$ is a 1-form on $U\cap V$ such that $dA=B_{U}-B_{V}$), one sees that the map $f$ in the sequence \eqref{red seqtk} can be expressed fully as the block matrix
$$
f=\left(
\begin{array}{cccc}
1 & h-1 & 0 & 0 
\\
0 & -h & 0 & 0
\\
0 & 0 & 1 & -e^{dA}
\end{array}
\right):\ZZ\oplus \ZZ \oplus \Omega^{\rm even}_{\rm cl}(U)\oplus 
\Omega^{\rm even}_{\rm cl}(V) \; \longrightarrow  \;
\ZZ\oplus \ZZ \oplus 
d\Omega^{\rm odd}(U\cap V)\;.
$$
The upper left block coming from the topological part is easily seen to be injective when $h\neq 0$. Thus, we need only understand the lower right block $(1, -e^{dA})$. Let 
$\omega=\omega_0+\omega_2$ and $v=v_0+v_2$ on $U$ and $V$ be a pair of closed forms on $U$ and $V$, respectively. \footnote{Recall that $U$ and $V$ have dimension 3, so there are no higher degree forms in $\omega$ and $v$.}
 Expanding $e^{dA}=1+dA$ on $U\cap V$, we see that the map $(1, -e^{dA})$ 
 is zero precisely if the equations
$$
\omega_0=v_0
\qquad \text{and} \qquad
\omega_2 = v_2-dA\wedge v_0
$$
are satisfied. Notice that if $H$ is not exact and $v_0\neq 0$, then no such solutions exist. Indeed, since $\omega$ and $v$ are closed, $v_0$ is constant and the second equation implies that
$
\omega_2/v_0-v_2/v_0=dA
$.
On the other hand, there are local potentials $B_{U}$ and $B_{V}$ for $H$ on $U$ and $V$, which satisfy $B_{U}-B_{V}=dA\;.$ Combining these equations, we have $\omega_2/v_0-B_{U}=v_2/v_0-B_{V}$ on $U\cap V$. By the sheaf condition for differential forms, there is a global form $\eta$, such that $\eta_{{}_{U}}=\omega_2/v_0-B_{U}$ and $\eta_{{}_{V}}=v_2/v_0-B_{V}$. But this means $d\eta=H$. Thus, if $H$ represents a nonzero cohomology class, then $v_0=\omega_0=0$ and $\omega_2=v_2$. By gluing, the pair $(\omega_2,v_2)$ are the restrictions of a global closed 2-form on $S^3$. We now have an exact sequence
\(\label{doubred seqtw}
\xymatrix{
K_{\RR/\ZZ}^{-2}(U; {h})\oplus K_{\RR/\ZZ}^{-2}(V; {h})\ar[r]^-{g} & 
K_{\RR/\ZZ}^{-2}(S^2; {h})\ar[r] & \widehat{K}^0(S^3; \hat{h})\ar[r] & \Omega^2_{\rm cl}(S^3)\to 0\;.
}
\)
Turning our attention to the map $g$, we again observe that since the twist $h$ trivializes on $S^2$, 
the Bockstein sequence for $K$-theory with coefficients, along with periodicity, gives an identification $K_{\RR/\ZZ}^{-2}(S^2; h)\cong \RR/\ZZ\oplus \RR/\ZZ$. As we have already observed, the map is induced by restrictions and the twist which in this case gives the map
$$
\left(\begin{array}{cc}
1 & h-1
\\
0 & -h
\end{array}\right):\RR/\ZZ\oplus \RR/\ZZ
\longrightarrow
 \RR/\ZZ\oplus \RR/\ZZ\;,
$$
where the matrix is to be understood as the map of $\ZZ$-modules which sends a pair $(\theta,\varphi)\in \RR/\ZZ$ to $(\theta (h-1)\varphi,-h\varphi)\in \RR/\ZZ$. This map is surjective and hence the sequence \eqref{doubred seqtw} reduces to 
\(
\xymatrix{
0\ar[r] & \widehat{K}^0(S^3; \hat{h})\ar[r] & \Omega^2_{\rm cl}(S^3)\ar[r] & 0
}\;,
\)
and we have an isomorphism 
$
\widehat{K}^0(S^3; \hat{h})\cong \Omega_{\rm cl}^2(S^3)
$.
\endofproof

\begin{proposition}
[Twisted differential $K^1$  of the 3-sphere]
\label{Prop-K1S3}
We have
$$
\widehat{K}^1(S^3; \hat{h})\cong \Omega^{\rm odd}(S^3)_{d_H-{\rm cl}}\;,
$$
where $d_H${--}${\rm cl}$ indicates closed forms with respect to the twisted 
de Rham differential $d_H$.
\end{proposition} 

\theproof
As in the case for degree zero above, we consider the portion of the Mayer-Vietoris sequence 
\(\label{K1red seqtk}
\xymatrix@=1.7em{
K_{\RR/\ZZ}^{-1}(S^2; h)\ar[r] & \widehat{K}^1(S^3; {\hat h})\ar[r] & 
\widehat{K}^1(U; {\hat h})\oplus \widehat{K}^1(V; {\hat h})\ar[r]^-{f} & \widehat{K}^1(U\cap V; {\hat h})\ar[r] & K^{2}(S^3; h)
}\;.
\)
By periodicity, this immediately reduces to 
\(\label{K1doubred seqtw}
\xymatrix{
0\ar[r] & \widehat{K}^1(S^3; {\hat h})\ar[r] & 
\widehat{K}^1(U; {\hat h})\oplus \widehat{K}^1(V; {\hat h})\ar[r]^-{f} & \widehat{K}^1(U\cap V; {\hat h})\ar[r] & 0
}\;.
\)
As before, we arrive at the following identifications
\begin{eqnarray*} 
\widehat{K}^1(U) &\cong& \RR/\ZZ\oplus \Omega^{\rm odd}_{\rm cl}(U)\;, 
\\
\widehat{K}^1(V) &\cong& \RR/\ZZ\oplus \Omega^{\rm odd}_{\rm cl}(V)\;, 
\\
\widehat{K}^1(U\cap V)&\cong& \RR/\ZZ\oplus \RR/\ZZ\oplus d\Omega^{\rm even}(U\cap V)\;, 
\end{eqnarray*}
and we identify the map $f$ with 
$$
f=\left(
\begin{array}{cccc}
1 & h-1 & 0 & 0 
\\
0 & -h & 0 & 0
\\
0 & 0 & 1 & -e^{dA}
\end{array}
\right):\RR/\ZZ\oplus \RR/\ZZ \oplus \Omega^{\rm odd}_{\rm cl}(U)\oplus 
\Omega^{\rm odd}_{\rm cl}(V) \; \longrightarrow  \;
\RR/\ZZ\oplus \RR/\ZZ \oplus 
d\Omega^{\rm even}(U\cap V)\;.
$$
In $f$, the upper left block matrix 
\(\label{up blk1}
\left(\begin{array}{cc}
1 & h-1
\\
0 & -h
\end{array}\right):\RR/\ZZ\oplus \RR/\ZZ
\longrightarrow
 \RR/\ZZ\oplus \RR/\ZZ
\)
is again understood as the map of $\ZZ$-modules which sends a pair $(\theta,\varphi)\in \RR/\ZZ$ to $(\theta (h-1)\varphi,-h\varphi)\in \RR/\ZZ$. Unlike in the degree zero case, however, 
now the kernel of the map \eqref{up blk1} is nonzero since we are taking $\RR/\ZZ$-coefficients. The requirement $-h\varphi\in \ZZ$ implies that (identifying $\RR/\ZZ\cong U(1)$) 
$\varphi$ is an $h$-root of unity. On the other hand, this combined with the second requirement $(1-h)\theta\varphi\in \ZZ$ implies that $\theta\in \ZZ$.
Hence we identify the kernel with the subgroup isomorphic to the image of the map
$$
\xymatrix{
\ZZ/h \; \ar@{^{(}->}^-{\times \frac{1}{h}}[r] & \;  \RR/\ZZ\oplus \RR/\ZZ\;,
}$$
sending $a\mapsto (1,a/h)$. Turning to the lower left block matrix $(1,-e^{dA})$ in $f$, we consider a pair of odd forms $\omega=\omega_1+\omega_3$ and $v=v_1+v_3$ on $U$ and $V$, respectively. Again,
 writing $e^{dA}=1+dA$, we identify the elements in the kernel as those forms satisfying the equations
$$
\omega_1=v_1
\qquad \text{and} \qquad
\omega_3 = v_3-dA\wedge v_1\;.
$$
Again we see that globally defined closed 3-forms give solutions (setting $v_1=0$). However we also have other solutions, even when $H$ is nontrivial in cohomology. The resulting forms are the twisted-closed
odd forms on $S^3$.
\endofproof

\medskip
Although the Mayer-Vietoris sequence provides a powerful tool for calculating twisted differential K-theory, the $\widehat{\rm AHSS}_{\hat h}$ can simplify some of these calculations drastically. We 
illustrate the power of  the spectral sequence in the following. 

\begin{example}
[The 3-sphere, revisited]
Let $\hat{h}:S^3\to \BB^2U(1)_{\nabla}$ be a twist. Recall from Proposition
 \ref{Prop-flat1st} that we identified the differential on the $E_3$-page as 
$d_3=\widehat{Sq}^3_{\ZZ}+\hat{h}\cup_{\rm DB}$. 
For the 3-sphere, the $U(1)$-cohomology is calculated from the exponential sequence as 
$$
H^2(S^3;U(1))\cong 0 \qquad \text{and} \qquad 
 H^3(S^3;U(1))\cong U(1)\;.
$$
Then for $\widehat{K}^0$, we see that all relevant differentials must vanish 
and the spectral sequence collapses at the $E_2$-page in Diagram 
\eqref{E2-Kth}. There is no extension 
problem in this case, and we have
$$
\widehat{K}^0(S^3; {\hat{h}})\cong (\Omega^{\rm even}(S^3),d_{H})_{\rm cl}\cong\Omega^2(S^3)_{\rm cl}\;.
$$
For $\widehat{K}^1$, we need to calculate the kernel of the map
$$
d_3:U(1)\longrightarrow U(1)\cong H^3(S^3;U(1))\;.
$$
For degree reasons, $\widehat{Sq}^3_{\ZZ}$ vanishes on $U(1)$ and we are reduced to finding the kernel of $\hat{h}$. To this end, recall that the formula for the Deligne-Beilinson cup product sends an element in $U(1)$, written as $e^{2\pi i \theta}$, to $e^{2\pi i h\theta}$, where $h$ is the integer representing the underlying topological twist. The kernel is thus easily identified with the $h$-roots of unity, which as an abelian group is isomorphic to $\ZZ/h\ZZ$. For degree reasons, there are no nontrivial differentials out of the term $(\Omega^*(S^3),d_{H})_{\rm cl}$. In this case, there is no extension problem and we recover the differential cohomology group, as in Proposition 
\ref{Prop-K1S3},
$$
\widehat{K}^1(S^3; {\hat{h}})\cong \ZZ/h\ZZ\oplus (\Omega^*(S^3),d_{H})_{\rm cl}\;.
$$
\end{example}

Generalizations to higher chromatic levels and higher spheres are studied 
in \cite{SY}, where twisted Morava K-theory of spheres are calculated. 
We plan to consider explicit refinements of those twisted theories elsewhere. 
We also plan to address many applications in a separate treatment, including 
the relation to classes in twisted Deligne cohomology \cite{GS4}, as well 
as implications for physics.


\end{document}